%% file: ksu.tex
\documentclass{amsart}
\usepackage{amssymb,epsfig, amsmath, amsfonts}

\usepackage{ifpdf}
\usepackage{color}
\usepackage{bm}
\usepackage{bbold} 
\usepackage{subfigure}
\usepackage{enumitem}
\usepackage{upgreek}
\usepackage{mathrsfs}
\usepackage[plain]{algorithm}
\usepackage{algorithmic}
\usepackage{framed}
\usepackage{afterpage}
\usepackage{graphicx}

\newtheorem{theorem}{Theorem}[section]
\newtheorem{lemma}[theorem]{Lemma}
\newtheorem{proposition}[theorem]{Proposition}

\theoremstyle{definition}
\newtheorem{definition}[theorem]{Definition}

\theoremstyle{remark}
\newtheorem{remark}[theorem]{Remark}

\numberwithin{equation}{section}

\floatstyle{ruled}
\newfloat{algo}{!ht}{loa}
\floatname{algo}{Algorithm}

\renewcommand{\algorithmicensure}{\emph{Parameters:}}

\usepackage{tikz}
\usepackage{pgfplots}
\usepackage{pgfplotstable}
\usetikzlibrary{calc}
\usetikzlibrary{positioning, fit}
\usetikzlibrary{decorations.pathreplacing}
\usetikzlibrary{shapes.symbols}
\usetikzlibrary{backgrounds}

\include{newcommands}

\title[A space-time adaptive wavelet method for time-periodic PDEs]{An efficient space-time adaptive wavelet Galerkin method for time-periodic parabolic partial differential equations}
\thanks{This work has partly been supported by the Deutsche Forschungsgemeinschaft within the Research Training Group (Graduiertenkolleg) GrK1100 \emph{Modellierung, Analyse und Simulation in der Wirtschaftsmathematik} at Ulm University.}
\date{\today}

\author{Sebastian Kestler}

\address{Sebastian Kestler,
Institute for Numerical Mathematics,
University of Ulm,
Helm\-holtz\-strasse 20,
D-89069 Ulm, Germany}

\email{sebastian.kestler@uni-ulm.de}

\author{Kristina Steih}
\address{
Kristina Steih,
Institute for Numerical Mathematics,
University of Ulm,
Helm\-holtz\-strasse 20,
D-89069 Ulm, Germany}
\email{kristina.steih@uni-ulm.de}

\author{Karsten Urban}
\address{
Karsten Urban,
Institute for Numerical Mathematics,
University of Ulm,
Helm\-holtz\-strasse 20,
D-89069 Ulm, Germany}
\email{karsten.urban@uni-ulm.de}

\makeatletter
\@namedef{subjclassname@2010}{%
  \textup{2010} Mathematics Subject Classification}
\makeatother



\begin{document}

\begin{abstract}
We {introduce} a multitree-based adaptive wavelet Galerkin algorithm {for} space-time discretized linear parabolic partial differential equations, focusing on time-periodic problems.
{It is shown that the} method converges with the best possible rate in linear complexity and can be applied for a wide range of wavelet bases.
We discuss the implementational challenges arising from the Petrov-Galerkin nature of the variational formulation and present numerical results for the heat and a convection-diffusion-reaction equation.
\end{abstract}

\maketitle

\section{Introduction}

In recent years, space-time {variational} approaches for linear parabolic partial differential equations (PDE) of the form
\[ u_{t} + \OpSp[u] = g(t) \qquad \text{ on } \Omega \subset \R^{n}, \text{ for } t \in [0,T], \]
have been {considered} in various contexts. These methods treat both temporal and spatial variables simultaneously, allowing e.g.\ for targeted adaptive refinement of the numerical discretization in the full space-time domain or efficient parallelization.
{On the other hand, this {in general} amounts to solving an $(n+1)$-dimensional problem.} 
This differs from standard time-stepping techniques for time-dependent PDEs, which are usually based on semi-discretization schemes: The \emph{vertical method of lines} requires the solution of a system of coupled ordinary differential equations that arise from a discretization in space. Within the \emph{horizontal method of lines} and the \emph{discontinuous Galerkin method}, the temporal variable is discretized first, leading to a (coupled) sequence of elliptic problems in the spatial domain. Such time-stepping schemes {have some} drawbacks: The sequential treatment of the time variable {often} {does not allow for}  parallelization in time.
Furthermore, adaptive schemes {typically} focus either on the spatial or on the temporal variable or are based on \emph{local} error estimators (e.g.\ \cite{Raasch:2007}), thus forfeiting optimality.
 Moreover, a posteriori error estimators -- needed e.g.\ for adaptive schemes or model reduction approaches -- are usually increasing functions in time, therefore losing efficiency over long time horizons. {These issues}  are amplified when considering time-periodic problems, i.e., when searching for solutions $u$ with $u(0) = u(T)$.
Such problems arise naturally in different physical, biological or chemical models, e.g.\ flows  around a rotor or propeller, biological models or chemical engineering \cite{VSPOpt,PeriodicalBioModel,Kawajiri:2006}. 
Standard numerical {methods for } periodic problems require either an additional fixed-point scheme (when using a temporal semi-discretization) or the solution of a system of boundary value problems (in case of the method of lines), both entailing non-negligible additional computational effort. 
In the present work, we will therefore consider a space-time variational formulation for time-periodic problems.

Space-time {variational} formulations for initial value problems {in particular} include space-time multigrid methods \cite{ST-Multigrid}, space-time sparse grids \cite{Andreev:2012,GriebelOeltz} or space-time wavelet collocation methods \cite{AlamKevlahanVasilyev}. Other space-time formulations based on special test bases or discontinuous Galerkin methods are e.g.\ \cite{MeidnerVexler,Urban2012203,UPFull}. 
These approaches exploit the space-time approach mainly for theoretical considerations and {allow}  the {use} of a time-stepping scheme -- thus effectively circumventing the main drawback of space-time methods, i.e., the additional dimension introduced by the temporal variable. {However, optimality has not yet been investigated in such a framework.}

Here, we follow the approach proposed by \cite{Schwab:2009}, where a space-time adaptive scheme using tensorized wavelet bases is proven to be optimal for initial value problems. In this setting, the partial differential equation is reformulated as an equivalent {\emph{non-symmetric}} bi-infinite matrix-vector problem of the form $\bB \bu = \bf$ and is numerically {approximated} by employing an adaptive wavelet Galerkin method (AWGM) to the corresponding normal equations.
{As opposed} to standard algorithms for time-periodic problems, the upshot of this approach is that time-periodic boundary conditions {can be} \emph{incorporated} into the underlying ansatz basis.

AWGMs {may be described as follows}, \cite{Cohen:2001,Gantumur:2007}:
{Consider} a bi-infinite linear system $\bC \bu = \bg$ in $\ell_2$ with a \emph{symmetric} positive definite {(s.p.d.)} stiffness matrix $\bC: \ell_2 \to \ell_2$, an infinite right-hand side $\bg \in \ell_2$ and a unique solution $\bu \in \ell_2$ which arises {e.g.}\ from the wavelet discretization of an elliptic operator problem.
In each iteration, these bi-infinite problems are {approximated} on a finite-dimensional index set $\bLambda_{k}$ steering the local refinement.
This means that a finite vector $\bu_{\bLambda_k}$ satisfying $\bC|_{\bLambda_k \times \bLambda_k} \bu_{\bLambda_k} = \bg |_{\bLambda_k}$ is computed. The {(infinite-dimensional) residual $\br_{\bLambda_{k}} :=\bg - \bC \bu_{\bLambda_{k}}$} is then approximated to serve as an error estimator and to identify an {update}, i.e.\  an ({usually} refined) index set $\bLambda_{k+1}$  {(typically $\bLambda_{k+1}$ corresponds to the significant coefficients of $\br_{\bLambda_{k}}$}). 

Space-time {variational} approaches to parabolic problems lead to \emph{non-symmetric} Petrov-Galerkin formulations and hence do not directly fit into this framework.
In particular, the bi-infinite associated stiffness matrix is no longer s.p.d.\ which is, however, a crucial ingredient for the convergence {analysis} of AWGMs.
Moreover,  the residual belonging to a test space which is \emph{not} identical to the trial space does not directly convey information {for an update of the \emph{trial} space}.
So, working with the normal equations, i.e., with the s.p.d.\ operator $\bC:=\bB^\top\bB$ and right-hand side $\bg = \bB^\top \bf$, is a natural approach for initial value problems {(e.g.\ \cite{Chegini:2011,Schwab:2009})} and, {as well we will show in this article,} also {for} {time-periodic} problems. 

The treatment of normal equations by adaptive wavelet methods has first been discussed in \cite{Cohen:2002}.
The main difficulty lies in the (approximate) evaluation of $\bB^\top \bB$ and $\bB^\top \bf$. 
%
There are several approaches in the literature that address these issues.
The techniques proposed in \cite{Cohen:2001} rely on so-called \emph{wavelet matrix compression schemes}, whereas \cite{Chegini:2011,Chegini:2012a}  use special wavelet constructions leading to truly sparse matrices $\bB$.
In the present work, we use another approach based on \emph{multitree-structured} index sets as introduced in \cite{Kestler:2012c,Kestler:2012d}.
This allows for the \emph{exact} application of  {$\bB$ for} wavelet discretizations of  linear differential operators with polynomial coefficients  within linear complexity when the underlying wavelet basis is of \emph{tensor product} type.
It basically consists of the so-called \emph{unidirectional principle} first introduced in sparse grid algorithms (e.g.\ \cite{Bungartz:2004,Zenger:1991}) where coordinate directions may be treated \emph{separately}.
The evaluation is then based on multitree-structured index sets that permit a tree structure when all but one coordinate directions are {frozen}. 

We stress that, {to the best of our knowledge}, \emph{no} quantitative results on the numerical solution of parabolic operator problems using the multitree concept within an AWGM  are available.
So far, only numerical results for \emph{elliptic} operator problems are presented in \cite{Kestler:2012c,Kestler:2012d}.
{Even though we focus on time-periodic problems}, we expect that our results can be extended to non-periodic settings {as well}.

The outline {of} this article is as follows.
In Section \ref{sec:problemsetting}, we introduce parabolic problems with periodic boundary conditions in time.
The derivation of the  equivalent $\ell_2$-problem by means of tensor product wavelet bases is explained  in Section \ref{sec:ell2problem}.
{Next, in Section \ref{sec:quasioptimal}, we define \emph{quasi-optimal} algorithms showing what can be expected in terms of convergence rates and computational work.} 
{Some} details on wavelet bases are then given in Section \ref{sec:waveletbases}.
In Section \ref{sec:awgm}, we explain AWGMs {for} elliptic  problems and extend it to our parabolic problem.
The realization and analysis of an efficient, multitree-based AWGM is then {presented in} Sections \ref{sec:mt_awgm} and \ref{sec:impl}.
We underline our theoretical findings by numerical experiments in Section \ref{sec:numerics}.

\section{{Time-periodic parabolic problems}}  \label{sec:problemsetting}

Let $\Domain  := \Omega_1 \times \cdots \times \Omega_n \subset \R^{n}$ {be} a product domain and $V$ {be} a {real} separable Hilbert space with dual $V'$ such that  $V \hookrightarrow H := L_{2}(\Domain) \hookrightarrow V'$ {is} a Gelfand triple. For $\OpSp \in \cL(V,V')$ and $g \in L_{2}(0,T;V')$ we consider the time-periodic equation
\begin{equation}\label{eq:problem_strongform}
		u_{t} + \OpSp[u(t)] = g(t) \, \text{ in } V' \text{ for a.e. }t \in [0,T],
		\qquad
		u(0) = u(T) \, \text{ in } H.  
\end{equation}
Denoting by $\eval{\cdot}{\cdot}_{V\times V'}$ the duality pairing on $V \times V'$, we assume that $t \mapsto \eval{v}{\OpSp[u]}$ is measurable on $[0,T]$ and {that} $\OpSp$ is uniformly coercive and bounded in time, i.e., there exist $0 < \alpha \leq \alpha(t)$, $\infty > \gamma \geq \gamma(t)$ such that for a.e. $t \in [0,T]$
\begin{equation} \label{eq:ellipticityA} 
	\eval{v}{\OpSp[w]}_{V\times V'} \leq \gamma \norm{w}_{V}\norm{v}_{V}, 
	\quad
	\eval{v}{\OpSp[v]}_{V\times V'} \geq \alpha \norm{v}_{V}^{2} 
	\hspace{0.5cm} \forall\, v,w \in V. 
\end{equation}
{Moreover}, we assume that the space $V$ is a Sobolev space of nonnegative order 
\begin{equation} \label{eq:intersection_V}
   V := \bigcap_{i=1}^n \bigotimes_{j=1}^n W_{ij}, \text{ where } W_{ij} :=
   \left\{
      \begin{array}{ll}
         L_2(\Omega_i), & i \neq j,\\
         V^{(i)}, & i = j,
      \end{array}
   \right.
\end{equation}
and, for a fixed $m \in \N$, $V^{(i)}$ is either $H^m(\Omega_i)$ or a subspace  incorporating essential boundary conditions.
{Note that several partial differential operators allow such a structure.}
As an example, think of $V= H^1_0(\bOmega)$ and $V^{(i)} = H^1_0(\Omega_i)$ (see \cite{Griebel:1995}).

\subsection{Space-time formulation}
We {derive a variational} formulation where the temporal {periodicity}  can be integrated into the function spaces and {is} therefore ultimately {incorporated} into the basis of a discrete approximation space. To this end, we derive a \emph{space-time variational formulation}: Defining 
\begin{align}\label{def:periodicsobspace}
   H^1_\per(0,T) &:= \{ v \in H^1(0,T): v(0) = v(T) \},
\end{align}
we consider the spaces $\cY := L_2(0,T; V )$ and $\cX := L_2 (0,T; V ) \cap H^1_{\per}(0,T; V')$, i.e.,
\begin{align}
   \cX &= \{ v \in L_2 (0,T; V ): v_{t} \in L_2 (0,T;V'),\ v(0) = v(T) \text{ in } H\}, \label{eq:cX} 
\end{align}
where $\cX$ is equipped with the norm
$\norm{v}^2_{\cX} := \norm{v}_{L_{2}(0,T;V)}^{2} + \norm{v_{t}}_{L_{2}(0,T;V')}^{2}$, $v \in \cX$.
Note that $v(0)$, $v(T)$ are well-defined due to $H^{1}(0,T) \subset {C([0,T])}$ and $\{ v \in L_2 (0,T; V ): v_{t} \in L_2 (0,T;V')\} \subset C(0,T;H)$, e.g.\ \cite{DautrayLions}. 
%
By integration of \eqref{eq:problem_strongform} over $[0,T]$, we obtain the problem:
\begin{align}\label{eq:problem_weakform}
 \text{Find } u \in \cX: \qquad b(u,v) = f(v) \qquad \forall\, v \in \cY,
\end{align}
with forms $b(\cdot, \cdot): \cX \times \cY \to \R$, $f(\cdot): \cY \to \R$ defined by, \cite[(5.6)-(5.7)]{Schwab:2009}
\begin{align}
  b(u,v) &:= \int_0^T [\eval{v(t)}{u_{t}(t)+\OpSp[u]}_{V\times V'} dt, \,
  \label{eq:bil_form_b}
  f(v) := \int_0^T \eval{v(t)}{g(t)}_{V\times V'}  dt.
\end{align}
We define the space-time operator $\OpST \in \cL(\cX,\cY')$ by $\eval{v}{\OpST[u]} := b(u,v)$ with  $\eval{\cdot}{\cdot} := \eval{\cdot}{\cdot}_{\cY \times \cY'}$, so that \eqref{eq:problem_weakform} is {a} variational formulation of the operator equation{:}
\begin{equation}\label{eq:problem_operatorform}
\text{Find } u \in \cX: \qquad \OpST[u] = f, \qquad f \in \cY'.
\end{equation}

\subsection{Well-posedness}
The well-posedness of a space-time formulation of (non-periodic) \emph{initial value} problems has been discussed in \cite{Schwab:2009}. In Appendix \ref{sec:appendix-wellposed}, we verify the Babu\v{s}ka-Aziz conditions:
\begin{enumerate}\renewcommand{\labelenumi}{(\roman{enumi})}
\item \emph{Continuity}: $\gammaB := \sup_{0 \neq u\in \cX} \sup_{0 \neq v \in \cY} \frac{b(u,v)}{\norm{u}_{\cX}\norm{v}_{\cY}} < \infty$.
\item \emph{Inf-sup condition:} $\betaB := \inf_{0 \neq u\in \cX} \sup_{0 \neq v \in \cY} \frac{b(u,v)}{\norm{u}_{\cX}\norm{v}_{\cY}} > 0$.
\item \emph{Surjectivity:} $\sup_{0 \neq u\in \cX} |b(u,v)| > 0$ for all $0 \neq v \in \cY$.
\end{enumerate}

\begin{proposition}\label{Prop:2.1}
Problem \eqref{eq:problem_weakform} is well-posed. In particular, $\OpST$ from \eqref{eq:problem_operatorform} is boundedly invertible with $\norm{\OpST} = \gammaB = \sqrt{2}\max\{1,\gamma\}$, $\norm{\OpST^{-1}} = \tfrac{1}{\betaB} = \frac{ \sqrt{2}\max\{1, \alpha^{-1}\} }{\alpha \min\{1, \gamma^{-2}\} }$.
\end{proposition}


\section{{Equivalent bi-infinite matrix-vector problem}} \label{sec:ell2problem}

We  {consider} the reformulation of \eqref{eq:problem_operatorform} as an \emph{equivalent} $\ell_2$-problem, i.e., a \emph{discrete} problem posed on the sequence space $\ell_2$. {This} was first introduced in \cite{Cohen:2001,Cohen:2002} for stationary problems and extended to parabolic problems in \cite{Schwab:2009}.

\subsection{Riesz bases} \label{subsec:riesz_bases}
We recall that for a separable Hilbert space $\cH$ of infinite dimension, a {dense} collection $\Upsilon := \{ \gamma_i : i \in \N \} \subset \cH$ is called a \emph{Riesz basis} for $\cH$ if there exist constants $\mathrm{c}, \mathrm{C} >0$ such that {for $v = \sum_{i=1}^\infty v_i \gamma_i$, it holds that}
\begin{equation} \label{eq:riesz_basis}
{  \mathrm{c}\| \bv \|^2_{\ell_2(\N)} \leq \| v \|^2_{\cH} \leq \mathrm{C} \| \bv \|^2_{\ell_2(\N)} \quad \forall \bv = (v_i)_{i \in \N} \in \ell_2(\N).   }
\end{equation}
The largest $\mathrm{c}$ and the smallest  $\mathrm{C}$ for which \eqref{eq:riesz_basis} holds, are referred to as \emph{lower} and \emph{upper Riesz constant} and are denoted by $\rmc{\Upsilon}{\cH}$ and 
$\rmC{\Upsilon}{\cH}$, respectively.\footnote{{Sometimes a different} definition of Riesz constants {is used}, namely $\rmc{\Upsilon}{\cH}$ and $\rmC{\Upsilon}{\cH}$ being the largest and the smallest constant such that $\rmc{\Upsilon}{\cH}\| \bv \|_{\ell_2(\N)} \leq \| v \|_{\cH} \leq \rmC{\Upsilon}{\cH} \| \bv \|_{\ell_2(\N)}$.}


\subsection{Wavelet discretization of the parabolic operator problem}

Let us now consider two \emph{different} Riesz bases
\begin{equation} \label{eq:Psis_cXcY}
   \hatbPsi^{\cX} := \big\{ \hatbpsi^{\cX}_{\blambda} : \blambda \in \hatbcJ \big\} \subset \cX, \qquad \checkbPsi^{\cY} := \big\{ \checkbpsi^{\cY}_{\blambda} : \blambda \in \checkbcJ \big\} \subset \cY,
\end{equation}
labeled w.r.t.\ two (possibly) different countable index sets $\hatbcJ$ {and} $\checkbcJ$.
More precisely, we consider a \emph{trial} basis $\hatbPsi^{\cX}$ for the ansatz space $\cX$ and a \emph{test} basis  $\checkbPsi^{\cY}$ for the test space $\cY$  with associated Riesz constants $\rmc{\cX}{\hatbPsi^{\cX}}$, $\rmC{\cX}{\hatbPsi^{\cX}}$ {and} $\rmc{\cY}{\checkbPsi^{\cY}}$, $\rmC{\cY}{\checkbPsi^{\cY}}$.
{It is important to note that $\hatbPsi^{\cX}$, $\checkbPsi^{\cY}$ arise from normalizing \emph{different} Riesz bases $\hatbPsi$, $\checkbPsi$, for $L_2((0,T)\times \bOmega)$ w.r.t.\ $\| \cdot \|_{\cX}$ and $\| \cdot \|_{\cY}$ ({which is also the reason for our notation}, see Section \ref{sec:waveletbases}).} 
%
{Then} there exists a \emph{unique} expansion $u = \bu^\top \hatbPsi^{\cX}$ of the solution $u$ of \eqref{eq:problem_operatorform} where we formally interpret both $\bu \in \ell_2(\hatbcJ)$ and $\hatbPsi^{\cX}$ as column vectors.
Now, the equivalent formulation of \eqref{eq:problem_operatorform}  reads {as} follows:
\begin{equation} \label{eq:equiv_ell2_prob}
  \text{Find }\bu \in \ell_2(\hatbcJ): \qquad \bB \bu = \bf,\qquad \bf \in \ell_2(\checkbcJ),
\end{equation}
where $\bB := \big[ \eval{\checkbpsi^{\cY}_{\blambda}}{\cB[\hatbpsi^{\cX}_{\bmu}]} \big]_{\blambda \in \checkbcJ, \bmu \in \hatbcJ} = \big[ b(\hatbpsi^{\cX}_{\bmu},  \checkbpsi^{\cY}_{\blambda} )\big]_{\blambda \in \checkbcJ, \bmu \in \hatbcJ}= \eval{\checkbPsi^{\cY}}{\cB[\hatbPsi^{\cX}]}$ is the \emph{bi-infinite stiffness matrix} and $\bf = \big[ \eval{\checkbpsi^{\cY}_{\blambda}}{f} \big]_{\blambda \in \checkbcJ}= \eval{\checkbPsi^{\cY}}{f}$ is the \emph{infinite right-hand side}. 
%
%
It is easy to see that {\eqref{eq:equiv_ell2_prob} is well-posed.
Since,} $\bf \in \ell_2(\checkbcJ)$ and $\cB \in \cL(\cX,\cY')$ is boundedly invertible, also $\bB \in \cL(\ell_2(\hatbcJ), \ell_2(\checkbcJ))$ is boundedly invertible.
In particular, with $\| \cdot \| := \| \cdot \|_{\ell_2 \to \ell_2}$ (compare \cite[(2.2) \& (2.3)]{Schwab:2009})
\begin{equation} \label{eq;bounds_bB_bBinv}
    \| \bB \|  \leq \| \cB \|_{\cX \to \cY'} \, \rmC{\cX}{\hatbPsi^{\cX}}^{\frac12} \, \rmC{\cY}{\checkbPsi^{\cY}}^{\frac12} , \quad  \| \bB^{-1} \| \leq \frac{ \| \cB^{-1} \|_{\cY' \to \cX} }{  \rmc{\cX}{\hatbPsi^{\cX}}^{\frac12} \, \rmc{\cY}{\checkbPsi^{\cY}}^{\frac12}}.
\end{equation}

\subsection{Further notations} \label{subsec:further_notations}
We need to restrict the bi-infinite matrices $\bB$ and $\bB^\top$ in both rows and columns. 
{For a pair $(\bLambda,\bcJ)$} with $\bcJ \in \{\hatbcJ,\checkbcJ \}$ and $\bLambda \subseteq \bcJ$, {set}
\begin{equation} \label{eq:extension_restriction_operator}
     \bE_{\bLambda}: \ell_2(\bLambda) \to \ell_2(\bcJ), \quad \text{ and } \quad \bR_{\bLambda} := \bE^\top_{\bLambda} : \ell_2(\bcJ) \to \ell_2(\bLambda),
\end{equation}
where $\bE_{\bLambda}$ is the trivial embedding, i.e., the extension of  $\bv_{\bLambda} \in \ell_2(\bLambda)$ by zeros to $\ell_2(\bcJ)$.
Consequently, its adjoint $\bR_{\bLambda}$ is the restriction of $\bv \in \ell_2(\bcJ)$ to $\bv|_{\bLambda} \in \ell_2(\bLambda)$. 
For $\hatbLambda \subseteq \hatbcJ$ and $\checkbLambda\subseteq \checkbcJ$, we define the following restriction of $\bB$ and $\bB^\top$:
\begin{equation} \label{eq:bB_finite}
 \bbB{\checkbLambda}{\hatbLambda} := \bR_{\checkbLambda} \, \bB \, \bE_{\hatbLambda}, 
 \quad 
 \bB_{\hatbLambda} := \bbB{\checkbcJ}{\hatbLambda}, 
 \quad
  \bbBtr{\hatbLambda}{\checkbLambda} \,:= \bR_{\hatbLambda} \, \bB^\top \, \bE_{\checkbLambda}, \quad {}_{\hatbLambda}\bB^\top := \bbBtr{\checkbcJ}{\hatbLambda}. 
\end{equation}
{{Finally},  $C \lesssim D$ means that $C$ can be bounded by a constant times $D$ and $C \gtrsim D$ is defined as $D \lesssim C$.}
In this setting, $C \eqsim D$ is defined as $C \lesssim D$ and $C \gtrsim D$.


\section{Quasi-optimal algorithms for bi-infinite matrix-vector problems} \label{sec:quasioptimal}

We may now focus on the \emph{approximate} solution of \eqref{eq:equiv_ell2_prob}.
To this end, we first discuss what can be \emph{expected} in terms of \emph{convergence rate} and \emph{complexity}.

\subsection{Best $\cN$-term approximation}
For a given number of degrees of freedom (d.o.f.) $\cN \in \N$, the best approximation $v_\cN$ of a function $v = \bv^\top \hatbPsi^{\cX} \in \cX$  in the basis $\hatbPsi^{\cX}$ with $\cN$ d.o.f.\ is a \emph{nonlinear, best} $\cN$-\emph{term approximation} (e.g.\ \cite{DeVore:1998}), i.e., $v_\cN=\arg\sigma_\cN(v)$, where the best $\cN$-term approximation error is defined as
$$
  \sigma_\cN(v) := \inf_{\{ \hatbLambda \in \hatbcJ : \# \hatbLambda = \cN \}} \inf_{ \{  v_\cN \in \Span\{ \hatbpsi^{\cX}_{\blambda } : \blambda \in \hatbLambda   \} \} } \| v - v_\cN \|_{\cX}.
$$
Since $\hatbPsi^{\cX}$ is a Riesz basis, it holds that $\| \bv - \bv_\cN \|_{\ell_2} \eqsim \sigma_\cN(v)$ where {$\bv_\cN$} {always denotes} an $\cN$-term approximation of {the vector} $\bv$ (i.e., the $\cN$ largest coefficients in modulus of $\bv$). 
As described in \cite{DeVore:1998}, it is meaningful to collect {all} vectors $\bv \in \ell_2(\hatbcJ)$ that permit an \emph{approximation rate} $s>0$ in the sense that $\| \bv - \bv_\cN \|_{\ell_2} \lesssim \cN^{-s}$ within the \emph{nonlinear approximation class} (compare \cite[(2)]{Stevenson:2009}):
\begin{equation} \label{eq:cAs}
   \cA^s \!:=\! \big\{  \bv \in \ell_2(\hatbcJ) \!:\! \| \bv \|_{\cA^s} \!:=\! \sup_{\eps >0} \eps \!\cdot\! \big[  \min\{ \cN \in \N_0 \!:\! \| \bv - \bv_\cN \|_{\ell_2(\hatbcJ)} \leq \eps  \}  \big]^s \!<\! \infty  \big\}.
\end{equation}
For a given $\bv \in \cA^s$ and $\eps>0$, the required number of degrees of freedom $\cN_\eps$ in order to obtain $\| \bv - \bv_{\cN_\eps} \|_{\ell_2} \leq \eps$ is bounded by
$
   \cN_\eps \leq \eps^{-1/s} \| \bv \|^{1/s}_{\cA^s}.
$ 
It is important to remark that this bound on $\cN_\eps$ is usually \emph{sharp} (see \cite[(3)]{Stevenson:2009}).

\subsection{Quasi-optimal algorithms} \label{subsec:quasioptimal}
Let us now assume that the solution $\bu\in\cA^s$ for some $s>0$ and that we want to approximate it with a target tolerance $\eps>0$.
The \emph{benchmark} is given by a best $\cN_\eps$-term approximation $\bu_{\cN_\eps}$ satisfying {$\sigma_{\cN_\eps}(\bu) = \| \bu - \bu_{\cN_\eps} \|_{\ell_2(\hatbcJ)} \leq \eps$} which is, however, in general \emph{not} computable.
So, we need to focus on the computation of a \emph{quasi-optimal} {approximation}  $\bu_\eps$: 
\begin{enumerate}[label=(O\arabic{*}), ref=(O\arabic{*})]
  \item \emph{Convergence rate}: $ \| \bu - \bu_\eps \|_{\ell_2(\hatbcJ)} \leq \eps$ and $\#  \supp \bu_\eps \lesssim \eps^{-1/s} \| \bv \|^{1/s}_{\cA^s}$. \label{item:conv_rate}
  \item \emph{Computational work}: The number of operations required for the computation of $\bu_\eps$ is of order $\cO( \eps^{-1/s} \| \bu \|_{\cA^s}^{1/s})$,  i.e., for any $\eps>0$, $\bu_\eps$ can be computed within \emph{linear} complexity, {recalling that $\cN_\eps\lesssim \eps^{-1/s} \| \bu \|_{\cA^s}^{1/s}$}. \label{item:comp_work}
\end{enumerate}

In order {to} realize \ref{item:comp_work}, we require the wavelet bases $\hatbPsi^{\cX}$ and $\checkbPsi^{\cY}$ to be of tensor product type which will be the topic of the next section.
	
\section{Tensor product wavelet bases} \label{sec:waveletbases}

{Recall} that $\cX$ and $\cY$ can be characterized as follows (see \cite{Griebel:1995}),
\begin{align}
   \cX  \eqsim \big[ L_2(0,T) \otimes V \big]  \cap \big[ H^1_\per(0,T) \otimes V' \big],  \quad  \cY  \eqsim  L_2(0,T) \otimes V. 
\end{align}
Furthermore, by the {definition} of $V$ in \eqref{eq:intersection_V}, the construction of $\hatbPsi^{\cX}$ and $\checkbPsi^{\cY}$ can be obtained by \emph{tensorization of univariate} wavelet bases.

\subsection{Uniformly local, piecewise polynomial wavelet bases}
Let us consider a univariate Sobolev space $\cH \in \{ H^1_\per(0,T)$, $V^{(1)}, \ldots,$ $V^{(n)}\}$ {with $V^{(i)} \subset L_2(\Omega_i)$} and  a {univariate} wavelet basis $\Psi$ for $L_2(\Omega)$ where {$\Omega\subset\R$} is either $(0,T)$ (if $\cH = H^1_\per(0,T)$) or $\Omega_i$ (if $\cH = V^{(i)}$, {recall $\Omega_i\subset\R$, i.e., w.l.o.g.\ $\Omega_i=(0,1)$}),
\begin{equation} \label{eq:Psi}
   \Psi = \bigcup_{j \in \N_0} \Psi_j = \{ \psi_\lambda : \lambda = (j,k) \in \cJ \}  \subset \cH,
\end{equation}
{as well as} $\Psi_j := \{ \psi_\lambda : \lambda \in \cJ_j \}$ and $\cJ_j := \{\lambda \in \cJ :  |\lambda| = j\}$.
Here, $|\lambda| := j \geq 0$ {denotes} the \emph{level} (steering the \emph{diameter} of the support of $\psi_{j,k}$ in the sense that $\diam (\supp \psi_{j,k}) \eqsim 2^{-j}$) and $k$ is a \emph{translation index} indicating the \emph{position} of $\supp \psi_{j,k}$.
{Note that the elements of $\Psi_0$ are not wavelets but scaling functions.} {For details on wavelets on the interval, we refer e.g.\ to \cite{Urban:WaveletBook}.} 
By the \emph{Wavelet Characterization Theorem} {\cite{DahmenActa}}, {if} the elements of $\Psi$ (and also those of the unique dual wavelet basis) are sufficiently smooth, the properly normalized collections
$\{ \psi_\lambda / \| \psi_\lambda\|_{\cH} : \lambda  \in \cJ \}$, $\{ \psi_\lambda / \| \psi_\lambda\|_{\cH'} : \lambda  \in \cJ \}$ are Riesz bases for {the Sobolev spaces} $\cH$ and $\cH'$, respectively.
Besides {that}, we shall  assume that $\Psi$ is a \emph{uniformly local, piecewise polynomial wavelet basis of order} $d \in \N$, i.e.:
{
\begin{enumerate}[label=(W\arabic{*}), ref=(W\arabic{*})]
   \item Local supports: $\diam(\supp \psi_\lambda) \eqsim 2^{-|\lambda|}$ for all $\lambda \in \cJ$. \label{eq:local_support}
   \item Level-wise finite number of overlaps: There exists $C \in \N$ \emph{independent} of $j \in \N_0$ such that
   $
   \sup_{\lambda \in \cJ_j} \#\{  \lambda' \in \cJ_{j} : |\supp \psi_\lambda \cap \supp \psi_{\lambda'}| >0 \} \leq  C
   $.
   \item Piecewise polynomials: For all $\lambda \in \cJ$, $\psi_\lambda$ is a piecewise polynomial of maximum degree $d-1$ and has $\wt d$ vanishing moments (except for scaling functions and {few} boundary adapted wavelets).
\end{enumerate}
}
Furthermore, we assume that the projection
$
Q_j [v]  :=  \sum_{\{\lambda \in \cJ : |\lambda|< j\}} v_\lambda \psi_\lambda
$
for $v= \sum_{\lambda \in \cJ} v_\lambda \psi_\lambda$
satisfies the following Jackson estimates 
$\| \Id - Q_j  \|_{H^{d}(\Omega) \cap \cH \to L_2(\Omega)} \lesssim 2^{-dj}$,    
$\| \Id - Q_j  \|_{H^{d}(\Omega) \cap \cH \to\cH} \lesssim 2^{-(d-m)j}$, 
$\| \Id - Q_j  \|_{H^{d}(\Omega) \cap \cH \to \cH'} \lesssim 2^{-(d+m)j}$, 
where $m=1$ if $\cH = H^1_{\per}(0,T)$.

\subsection{Temporal discretization}
In order to {ensure} the periodic boundary conditions in time (see \eqref{eq:problem_strongform}) in $\cX$ {we need} a ({univariate}) \emph{periodic} wavelet basis 
\begin{equation} \label{eq:Thetaper}
   \Theta^{\per} := \big\{  \theta^{\per}_\lambda : \lambda \in \cJ^{\per}_t \big\} \subset H^1_\per(0,T)
\end{equation}
being a uniformly local, piecewise polynomial wavelet basis of order $d_t \in \N$ {(the index $t$ stands for `time')}  for $L_2(0,T)$ with associated Riesz constants $\rmc{L_2}{\Theta^\per}$, $\rmC{L_2}{\Theta^\per}$.
We assume that the elements of $\Theta^{\per}$ are sufficiently smooth so that the \emph{properly normalized} collection $\big\{  \theta^{\per}_\lambda/ \| \theta^{\per}_\lambda \|_{H^1} : \lambda \in \cJ^{\per}_t \big\}$ is a Riesz basis for $H^1_\per(0,T)$ with constants $\rmc{H^1_\per}{\Theta^\per}$, $\rmC{H^1_\per}{\Theta^\per}$. {Recall that the construction of \emph{periodic} wavelet bases is particularly easy, \cite{Urban:WaveletBook}.} 
For {the} temporal part of the \emph{test} space $\cY$ {(involving also non-periodic functions)}, we consider a uniformly local, piecewise polynomial wavelet basis for $L_2(0,T)$,
\begin{equation} \label{eq:Theta}
\Theta := \big\{  \vartheta_\lambda : \lambda \in \cJ_t \big\},
\end{equation}
with Riesz constants $\rmc{L_2}{\Theta}$, $\rmC{L_2}{\Theta}$ {and wavelets being  \emph{not} necessarily periodic}.

\subsection{Spatial discretization}
For the spatial discretization, we use the fact that $\Domain = \Omega_1 \times \cdots \times \Omega_n$ is a product domain.
Here, we shall use that $V$ is the (intersection of) tensor products of univariate Sobolev spaces (see \eqref{eq:intersection_V}) with $L_2(\Domain) \subseteq V$ and $L_2(\Domain) \eqsim L_2(\Omega_1) \otimes \cdots \otimes L_2(\Omega_n)$ (see, e.g., \cite{Griebel:1995}).
We assume that for $i \in \{1,\ldots,n\}$ we are given univariate uniformly local, piecewise polynomial wavelet bases of order $d_x \in \N$ {(the index $x$ indicating the spatial variable)} for $L_2(\Omega_i)$,
$\Sigma^{(i)} := \{ \sigma^{(i)}_\lambda : \lambda \in \cJ^{(i)} \} \subset V^{(i)}$. 
We require that these functions are sufficiently smooth so that  $\{ \sigma^{(i)}_\lambda/ \| \sigma^{(i)}_\lambda \|_{V^{(i)}} : \lambda \in \cJ^{(i)} \}$, $\{ \sigma^{(i)}_\lambda/ \| \sigma^{(i)}_\lambda \|_{V^{(i)}{}'}: \lambda \in \cJ^{(i)} \}$ are Riesz bases for $V^{(i)}$, $V^{(i)}{}'$ with  constants $\rmc{V^{(i)}}{\Sigma^{(i)}}$, $\rmC{V^{(i)}}{\Sigma^{(i)}}$ {and} $\rmc{V^{(i)}{}'}{\Sigma^{(i)}}$, $\rmC{V^{(i)}{}'}{\Sigma^{(i)}}$. 
Now, 
\begin{equation} \label{eq:bSigma}
   \bSigma := \big\{ \bsigma_{\blambda} : \blambda \in \bcJ_{\hspace{-0.5ex} x} \big\} := \Sigma^{(1)} \otimes \cdots \otimes \Sigma^{(n)}
\end{equation}
is a Riesz basis for $L_2(\Domain)$ where $\bsigma_{\blambda} := \sigma^{(1)}_{\lambda_1} \otimes \cdots \otimes \sigma^{(n)}_{\lambda_n}$ is a \emph{tensor product wavelet} and $ \bcJ_{\hspace{-0.5ex} x} := \cJ^{(1)} \times \cdots \times \cJ^{(n)}$, \cite[Lemma 3.1.7]{Dijkema:Diss}.
Moreover, 
\begin{equation} \label{eq:bSigma_normalized}
  \bSigma^{V} :=  \big\{ \bsigma_{\blambda} / \| \bsigma_{\blambda} \|_V : \blambda \in \bcJ_{\hspace{-0.5ex} x} \big\}, \quad   \bSigma^{V'} :=   \big\{ \bsigma_{\blambda} / \| \bsigma_{\blambda} \|_{V'} : \blambda \in \bcJ_{\hspace{-0.5ex} x} \big\}
\end{equation}
are Riesz bases for $V$, $V'$, \cite[Lemma 3.1.8]{Dijkema:Diss}.
The associated Riesz constants will be denoted by $\rmc{V}{\bSigma}$, $\rmC{V}{\bSigma}$, $\rmc{V'}{\bSigma}$ and $\rmC{V'}{\bSigma}$.

\subsection{Space-time discretization}\label{Sec:STDiscr}
We are now in the position to define the Riesz wavelet bases $\hatbPsi^{\cX}$ and $\checkbPsi^{\cY}$ from \eqref{eq:Psis_cXcY}.
With  $L_2(0,T; L_2(\Domain)) \eqsim L_2(0,T) \otimes L_2(\Domain)$,
\begin{align}
 &  \hatbPsi := \big\{  \hatbpsi_{\blambda} :=  \theta^{\per}_{\lambda_t}  \otimes \bsigma_{\blambda_x}: \blambda := (\lambda_t,\blambda_x) \in \hatbcJ := \cJ^{\per}_t \times \bcJx \big\} =    \Theta^{\per} \otimes \bSigma, \label{eq:hatPsiXL2} \\
&     \checkbPsi := \big\{  \checkbpsi_{\blambda} :=  \vartheta_{\lambda_t}  \otimes \bsigma_{\blambda_x}: \blambda := (\lambda_t,\blambda_x) \in \checkbcJ := \cJ_t \times \bcJx \big\} =    \Theta \otimes \bSigma, \label{eq:hatPsiYL2}
\end{align}
are both Riesz bases for $L_2(0,T; L_2(\Domain))$.
At this point, we only need to normalize the above Riesz bases appropriately (see \cite[Propositions 1 \& 2]{Griebel:1995}) so that
\begin{align}
&    \hatbPsi^{\cX} := \big\{   \hatbpsi_{\blambda} /   \| \hatbpsi_{\blambda} \|_{\cX} : \blambda \in \hatbcJ  \big\}  = \bD^{\cX} \hatbPsi, \quad \bD^{\cX} := \diag\big[ \big( \| \hatbpsi_{\blambda} \|^{-1}_{\cX} \big)_{\blambda \in \hatbcJ}\big],   \label{eq:hatbPsiX}\\
&    \checkbPsi^{\cY} := \big\{   \checkbpsi_{\blambda} /   \| \checkbpsi_{\blambda} \|_{\cY} : \blambda \in \checkbcJ  \big\}  = \bD^{\cY} \checkbPsi, \quad\; \bD^{\cY} := \diag\big[ \big( \| \checkbpsi_{\blambda}\|^{-1}_{\cY}  \big)_{\blambda \in \checkbcJ}\big],  \label{eq:checkbPsiY}
\end{align}
are Riesz bases for $\cX$, respectively $\cY$ (compare \cite[Section 6]{Schwab:2009}).

\begin{remark}
We shall denote a tensor product wavelet basis $\bPsi \in \{ \hatbPsi, \checkbPsi\}$ as follows:
\begin{align*}
    \bPsi &= \Psi^{(0)} \otimes \Psi^{(1)} \otimes \cdots \otimes \Psi^{(n)} = \big\{ \bpsi_{\blambda} := \psi^{(0)}_{\lambda_0} \otimes \psi^{(1)}_{\lambda_1} \otimes \cdots \otimes \psi^{(n)}_{\lambda_n}: \blambda \in \bcJ \big\},
\end{align*}
where $\blambda = (\lambda_0, \lambda_1,\ldots,\lambda_n)$ and $\bcJ := \cJ^{(0)} \times \cJ^{(1)} \times \cdots \times \cJ^{(n)}$.
In this setting, it is clear that $\Psi^{(0)} \in \{ \Theta^{\per}, \Theta\}$, $\cJ^{(0)} \in \{ \cJ^{\per}_t, \cJ_t \}$ and $\Psi^{(i)} = \Sigma^{(i)}$ for $i \in \{1,\ldots,n\}$.
\end{remark}

\subsection{Riesz constants for test and trial bases} \label{sec:condition_numbers}
For the implementation of an AWGM, we need \emph{estimates} for the Riesz constants $\rmc{\cX}{\hatbPsi}$, $\rmC{\cX}{\hatbPsi}$, $\rmc{\cY}{\checkbPsi}$, $\rmC{\cY}{\checkbPsi}$ in \eqref{eq;bounds_bB_bBinv}.
Again, we use that $\cX$ and $\cY$ are (intersections of) tensor products of Hilbert spaces.
As in \cite[\S6]{Schwab:2009}, we have the following estimates for $\hatbPsi^{\cX}$ and $\checkbPsi^{\cY}$
\begin{align}
&    \, \rmc{\cX}{\hatbPsi} \geq \min\big\{ \rmc{L_2}{\Theta^{\per}} \cdot \rmc{V}{\bSigma}, \, \rmc{H^1_\per}{\Theta^{\per}} \cdot \rmc{V'}{\bSigma}   \big\},  \label{eq:rmc_cX} \\
&     \rmC{\cX}{\hatbPsi} \leq \min\big\{ \rmC{L_2}{\Theta^{\per}} \cdot \rmC{V}{\bSigma}, \, \rmC{H^1_\per}{\Theta^{\per}} \cdot \rmC{V'}{\bSigma}   \big\}, \label{eq:rmC_cX} \\
&     \rmc{\cY}{\checkbPsi} \geq  \rmc{L_2}{\Theta} \cdot \rmc{V}{\bSigma}, \quad   \rmC{\cY}{\checkbPsi} \leq  \rmC{L_2}{\Theta} \cdot \rmC{V}{\bSigma}.  \label{eq:rmcC_cY} 
\end{align}

The Riesz constants $\rmc{V}{\bSigma}$, $\rmC{V}{\bSigma}$ can also be bounded by those of the 1D bases  $\Sigma^{(i)}$, $i \in \{1,\ldots,n\}$.
Using \eqref{eq:intersection_V},  it can be shown as in \cite[\S2]{Dijkema:2009}{, that}
\begin{align}
 &   \,\rmc{V}{\bSigma} \geq \min_{m \in \{1,\ldots,n\}} \min \Big\{ \rmc{L_2}{\Sigma^{(m)}} , \rmc{V^{(m)}}{\Sigma^{(m)}} \Big\} \prod_{k\neq m} \rmc{L_2}{\Sigma^{(k)}} , \label{eq:rmc_V} \\
 &  \rmC{V}{\bSigma} \leq \max_{m \in \{1,\ldots,n\}} \max \Big\{ \rmC{L_2}{\Sigma^{(m)}} , \rmC{V^{(m)}}{\Sigma^{(m)}}\Big\}  \prod_{k\neq m} \rmC{L_2}{\Sigma^{(k)}}. \label{eq:rmC_V}
\end{align}
Unfortunately, the same approach does not apply to the {(dual)} Riesz constants $\rmc{V'}{\bSigma}$, $\rmC{V'}{\bSigma}$ of $\bSigma^{V'}$ in \eqref{eq:bSigma_normalized}.
However, one may consider $\wt \bSigma^{V}$ being the unique Riesz basis for $V$ that is \emph{dual} to $\bSigma^{V'}$, i.e., $\eval{\wt \bSigma^{V}}{\bSigma^{V'}}_{V \times V'} = \Id$.
{Denoting by $\rmc{V}{\wt \bSigma}$, $\rmC{V}{\wt \bSigma}$ the associated Riesz constants, it can be shown that $\rmC{V}{\wt \bSigma}^{-1} \leq \rmc{V'}{\bSigma}$ and $\rmC{V'}{\bSigma} \leq \rmc{V}{\wt \bSigma}^{-1}$. 
Observe that for computing bounds for $\rmc{V}{\wt \bSigma}$, $\rmC{V}{\wt \bSigma}$, we may proceed as for bounding $\rmc{V}{\bSigma}$, $\rmC{V}{\bSigma}$.} 
We conclude that for the computation of the bounds in \eqref{eq:rmc_cX}, \eqref{eq:rmC_cX} and \eqref{eq:rmcC_cY}, it is sufficient to compute bounds for \emph{univariate} Riesz constants which can be easily approximated  (e.g.\ \cite[\S2]{Dijkema:Diss}).

\begin{remark}
{Recalling} the construction of wavelets, {note} that the \emph{numerical approximation} of $\rmc{V}{\wt \bSigma}$, $\rmC{V}{\wt \bSigma}$ {may} be difficult since the the dual basis $\wt \bSigma^{V}$ (and their derivatives) may not be available in a closed form. If sharp bounds are needed, one may use an $L_2(\Domain)$-orthonormal basis $\bSigma$ so that $\wt \bSigma^{V} = \bSigma^V$, {e.g.\ multiwavelets}.
\end{remark}

\subsection{Best approximation rates}
We need to know for which values of $s$ the solution $\bu$ of \eqref{eq:equiv_ell2_prob} is in $\cA^s$. 
More precisely, for a fixed trial basis $\hatbPsi^{\cX}$, the question is what is the \emph{largest} value $\smax$ of $s$ for which $\bu \in \cA^s$  can be expected and that cannot be increased by imposing higher smoothness conditions on $u$ (excluding special cases  where $\bu$ is (close to) a finite vector).
This value $\smax$ is referred to as \emph{best possible approximation rate}. 
For our setting, we may apply the results from \cite[\S7.2]{Schwab:2009}.
With $u = \bu^\top \hatbPsi^{\cX} \in \cX \cap H^{d_t}(0,T) \otimes \scrH^{d_x}(\bOmega)$ and the Sobolev space
$$
   \scrH^{d_x}(\bOmega) := \bigcap_{i=1}^n \bigotimes_{j=1}^n Z_{ij}, \text{ where } Z_{ij} :=
   \left\{
      \begin{array}{ll}
          L_2(\Omega_i), & i \neq j,\\
          {H^{d_x}(\Omega_i)}, & i = j, 
      \end{array}
   \right.
$$
of dominating mixed derivatives, the best possible rate is given by
\begin{equation} \label{eq:smax}
   \smax = \min\{ d_t -1 , d_x -m  \}.
\end{equation}
{We recall that}  $d_t$ denotes the polynomial order of $\Theta^{\per}$ and $d_x$ those of $\Sigma^{(1)}, \ldots, \Sigma^{(n)}$.
This rate does \emph{not} depend on the spatial dimension $n$.
Moreover, we remark that $u\in\scrH^{d_x}(\bOmega)$ is sufficient but not necessary {for obtaining the above rate}.
In fact, the Sobolev space $H^{d_t}(0,T) \otimes \scrH^{d_x}(\bOmega)$ can  be replaced by a (weaker) Besov space of dominating mixed derivatives, \cite{Nitsche:2006, Sickel:2009}. 
Note that the order of the wavelet bases for the \emph{test} space $\cY$ does \emph{not} enter the best approximation rate. 

\section{{Adaptive wavelet Galerkin methods}} \label{sec:awgm}

An infinite $\ell_2$-problem \eqref{eq:equiv_ell2_prob} arising from a wavelet discretization of  \eqref{eq:problem_operatorform} can be solved by an AWGM, e.g.\ \cite{Cohen:2001,Gantumur:2007}. We now first present the main idea of an AWGM for the solution of an (for convenience) \emph{elliptic} operator problem. 
Secondly, we highlight the additional {challenges} related to \emph{parabolic} problems and indicate a possible way-out using \emph{normal equations}.

\subsection{Elliptic operator problems} \label{subsec:awgm_elliptic}
\emph{Solely} for explanation purposes, we consider \emph{elliptic} operator problems of the following type.
For a linear, self-adjoint operator $\cC \in \cL(\cX,\cX')$ induced by a continuous and coercive bilinear form  (i.e., $\eval{v}{\cC[w]}_{\cX \times \cX'} \lesssim \| v \|_{\cX} \| w \|_{\cX}$, $\eval{v}{\cC[v]}_{\cX \times \cX'} \gtrsim \| v \|_{\cX}^{2}$ for all $v,w \in \cX$),
we consider:
\begin{equation} \label{eq:cUug}
  \text{Find } u \in \cX: \qquad   \cC[u] = g, \qquad g \in \cX'.
\end{equation}

Analogously to \eqref{eq:equiv_ell2_prob},  the equivalent $\ell_2$-problem to this problem reads:
\begin{equation} \label{eq:bCbubg}
  \text{Find } \bu \in \ell_2(\hatbcJ): \qquad  \bC \bu = \bg, \qquad \bg \in \ell_2(\hatbcJ),
\end{equation}
where $\bC = \eval{\hatbPsi^{\cX}}{\cC[\hatbPsi^{\cX}]}_{\cX \times \cX'}$ and $\bg= \eval{\hatbPsi^{\cX}}{g}_{\cX \times \cX'}$ with $\hatbPsi^{\cX}$ from \eqref{eq:hatbPsiX}.
%
{In the elliptic case}, i.e., $\cX = \cY$ and may use $\hatbPsi^{\cX}$ as trial \emph{and} test basis.
Furthermore, $\bC$ is s.p.d.\ and $\triplenorm \cdot \triplenorm^2 := \eval{\,\cdot\,}{\bC \,\cdot \,}_{\ell_2(\hatbcJ) \times \ell_2(\hatbcJ)}$ defines an equivalent norm,  \cite[p.\ 565]{Stevenson:2009}
\begin{equation} \label{eq:norm_bc}
   \| \bC^{-1} \|^{-\frac12}\| \bv \|_{\ell_2}  \leq \triplenorm \bv \triplenorm \leq \| \bC \|^{\frac12} \| \bv \|_{\ell_2}, \quad \forall \bv \in \ell_2(\hatbcJ).
\end{equation}
The idea of an AWGM for \eqref{eq:bCbubg} is outlined in an (idealized) Algorithm \ref{alg:exact_awgm}, \cite[p.\ 567]{Stevenson:2009}. 
Within this algorithm, we  make some non-realistic assumptions, which will be discussed below.
Abandoning these assumptions will then give rise to the realizable AWGM variants introduced in later sections. 
Starting from an initial index set $\hatbLambda_1\subset \hatbcJ$, a sequence of \emph{nested} finite index sets $(\hatbLambda_k)_k$ is computed.
On each {such} $\hatbLambda_k$, a Galerkin problem is solved that yields the (finite) vector $\bu_{\hatbLambda_k}$.
Due to the Riesz basis property, {it holds that} (see also \eqref{eq:riesz_basis})
$$
{\rmc{\cX}{\hatbPsi}^{\frac12}   \|   \bu -  \bu_{\hatbLambda_k} \|_{\ell_2(\hatbcJ)} \leq   \| u - \bu_{\hatbLambda_k}^\top \hatbPsi^{\cX}  \|_{\cX} \leq \rmC{\cX}{\hatbPsi}^{\frac12}   \|   \bu -  \bu_{\hatbLambda_k} \|_{\ell_2(\hatbcJ)}.}
$$
Given $\bu_{\hatbLambda_k}$, the computation of the next $\hatbLambda_{k+1}$ is based on the \emph{infinitely} {supported} \emph{residual} $\bg - \bC \bu_{\hatbLambda_k} \in \ell_2(\hatbcJ)$  {and the} \emph{error estimator} $\| \bg - \bC \bu_{\hatbLambda_k} \|_{\ell_2(\hatbcJ)}$
which satisfies:
\begin{equation} \label{eq:estimate_elliptic_problem}
    \| \bC \|^{-1} \| \bg - \bC \bu_{\hatbLambda_k} \|_{\ell_2(\hatbcJ)} \leq \|   \bu -  \bu_{\hatbLambda_k} \|_{\ell_2(\hatbcJ)}  \leq  \| \bC^{-1} \| \| \bg - \bC \bu_{\hatbLambda_k} \|_{\ell_2(\hatbcJ)}.
\end{equation}
This also explains the {stopping} criterion  in line 4 of Algorithm \ref{alg:exact_awgm}.
{Consequently}, indices corresponding to the largest entries in the residual are \emph{added} to $\bLambda_k$.
This so-called \emph{bulk-chasing} process is steered by the parameter $\delta$.

\renewcommand{\algorithmicensure}{\emph{Parameter:}}
\begin{algo}[!ht]
\caption{[$\bu_\eps$] = \textbf{IDEALIZED--AWGM}[$\eps$, $\bLambda_1$]
	\label{alg:exact_awgm}}
\algsetup{indent=2em}
\begin{algorithmic}[1]
\REQUIRE Target tolerance $\eps$ and an index set $\hatbLambda_1 \neq \emptyset$.
\ENSURE $\delta \in (0,\kappa(\bC)^{-\frac12})$.
\FOR{$k=1,2,\ldots$}
\STATE \emph{{Solve the} Galerkin problem:}
\begin{equation} \label{eq:galerkin_system_elliptic}
 \text{Find } \bu_{\hatbLambda_k} \in \ell_2(\hatbLambda_k): \qquad {}_{\hatbLambda_k} \bC_{\hatbLambda_k} \bu_{\hatbLambda_k} = \bg_{\hatbLambda_{k}}, \quad \bg_{\hatbLambda_k} := \bR_{\hatbLambda_k} \bg \in \ell_2(\hatbLambda_k).
\end{equation}
\STATE \emph{Residual computation:} Compute $\bg - \bC \bu_{\hatbLambda_k}$ and $\nu_k := \| \bg - \bC \bu_{\hatbLambda_k} \|_{\ell_2}$.
\STATE \textbf{if} $\nu_k \leq \eps / \| \bC^{-1} \|$ \textbf{then} \textbf{return} $\bu_\eps := \bu_{\hatbLambda_k}$.
\STATE \emph{Bulk chasing criterion:} Find smallest index set $\hatbLambda_{k+1} \supset \hatbLambda_k$ such that
\begin{equation} \label{eq:bulkchasing_system_elliptic}
  \| \bR_{\hatbLambda_{k+1}}(\bg - \bC \bu_{\hatbLambda_k}) \|_{\ell_2(\hatbLambda_{k+1})} \geq \delta \| \bg - \bC \bu_{\hatbLambda_k} \|_{\ell_2(\hatbcJ)}.
\end{equation}
\ENDFOR
\end{algorithmic}
\end{algo}
\renewcommand{\algorithmicensure}{\emph{Parameters:}}


\begin{proposition}[{\cite[Proposition 4.1]{Stevenson:2009}}] \label{prop:conv_idealized_awgm}
The {iterates $\bu_{\hatbLambda_k}$} produced by Algorithm \ref{alg:exact_awgm} satisfy
$\triplenorm \bu - \bu_{\hatbLambda_k} \triplenorm \leq [1 - \delta^2 \kappa(\bC)^{-1}]^{k/2} \triplenorm \bu \triplenorm$.
{For the output $\bu_\eps$ it holds } $\| \bu - \bu_\eps \|_{\ell_2(\hatbcJ)} \leq \eps$.
If $\bu \in \cA^s$ for some $s>0$, it also holds for $\cN_k := \# \hatbLambda_k$ that
\begin{equation}
     \| \bu - \bu_{\hatbLambda_k} \|_{\ell_2(\hatbcJ)} \lesssim \| \bu \|^{1/s}_{\cA^s} \cN_k^{-s} , \quad  \# \supp \bu_\eps \lesssim \eps^{-1/s} \| \bu \|^{1/s}_{\cA^s}.
\end{equation}
\end{proposition}

\begin{remark}
Algorithm \ref{alg:exact_awgm} cannot be implemented as the residual cannot be computed exactly in general.
\emph{Implementable} versions are given in \cite{Cohen:2001,Gantumur:2007}.
The algorithm in \cite{Cohen:2001} requires an additional \emph{thresholding} and thus can be expected to be  less efficient than \cite{Gantumur:2007}. 
The adaptive wavelet method in \cite{Cohen:2002} relies on an \emph{inexact} Richardson iteration that is applied \emph{directly} to \eqref{eq:bCbubg} without Galerkin projection.
However, as shown in \cite{Gantumur:2007}, also this scheme can be expected to be less efficient than \cite{Gantumur:2007}.
Thus, we shall  focus on \cite{Gantumur:2007} here.
\end{remark}

\subsection{Parabolic problems} \label{subsec:extension_idealawgm_parabolic}
One may try to analyze \textbf{IDEALIZED--AWGM} for $\bB \bu = \bf$ in \eqref{eq:equiv_ell2_prob}. However, the generalization of the idealized scheme to \eqref{eq:equiv_ell2_prob} is not trivial: 
\emph{(1) Symmetry and positive definiteness.}
Recall that $\bB$ from 
\eqref{eq:equiv_ell2_prob} is \emph{not} s.p.d., so that $\eval{\,\cdot\,}{\bB\,\cdot\,}_{\ell_2(\hatbcJ) \times \ell_2(\checkbcJ)}$ is \emph{not} an equivalent norm on $\ell_2(\hatbcJ)$.
However, the availability of an equivalent \emph{energy norm} as in \eqref{eq:norm_bc} is crucial for  the convergence analysis of Algorithm \ref{alg:exact_awgm} (see \cite[Proposition 4.1]{Stevenson:2009}). 
\emph{(2) Bulk chasing and residual computation.}
It is not clear how to construct $\hatbLambda_{k+1}$ from $\hatbLambda_k$.
In analogy to \eqref{eq:estimate_elliptic_problem}, the residual $\bf - \bB \bu_{\hatbLambda_k} \in \ell_2(\checkbcJ)$  with error estimator $\|  \bf - \bB \bu_{\hatbLambda_k} \|_{\ell_2(\checkbcJ)}$ satisfies
\begin{equation} \label{eq:estimate_nonelliptic_problem}
  \| \bB \|^{-1} \| \bf - \bB \bu_{\hatbLambda_k} \|_{\ell_2(\checkbcJ)} \leq  \|   \bu -  \bu_{\hatbLambda_k} \|_{\ell_2(\hatbcJ)} \leq   \| \bB^{-1} \| \| \bf - \bB \bu_{\hatbLambda_k} \|_{\ell_2(\checkbcJ)}.
\end{equation}
But the residual is an element of $\ell_2(\checkbcJ)$, $\checkbcJ \neq \hatbcJ$. Thus, we \emph{cannot} compute $\hatbLambda_k$ by selecting some contributions from the residual as in \eqref{eq:bulkchasing_system_elliptic}. 
\emph{(3) Petrov-Galerkin problems.} 
Since $\hatbPsi^{\cX} \neq \checkbPsi^{\cY}$, the (well-posed) Galerkin problem in line 5 of Algorithm \ref{alg:exact_awgm} here becomes a \emph{Petrov}-Galerkin problem. Hence the \emph{uniform well-posedness} of the finite-dimensional problems is no longer inherited from the infinite dimensional problem \eqref{eq:equiv_ell2_prob} and has to be ensured for all  $\hatbLambda_{k}$.
\medskip


Hence, we focus on the associated \emph{normal equations}, as proposed in \cite{Cohen:2002}:
\begin{equation} \label{eq:normal_equations}
  \text{Find } \bu \in \ell_2(\hatbcJ): \qquad \bB^\top \bB \bu = \bB^\top \bf, \qquad \bB^\top\bf \in \ell_2(\hatbcJ).
\end{equation}
Since $\bB$ is boundedly invertible, the unique solution of \eqref{eq:normal_equations} is also the unique solution of \eqref{eq:equiv_ell2_prob} (see \cite[Thm.\ 7.1]{Cohen:2002}).
Indeed, \eqref{eq:normal_equations} are the (infinite) \emph{normal equations} associated to the least squares problem (compare \cite[\S7]{Cohen:2002}) of finding 
$\bu \in \ell_2(\hatbcJ)$ such that
$\bu = \argmin_{\bv \in \ell_2(\hatbcJ)} \| \bB \bv - \bf \|_{\ell_2(\checkbcJ)}^2$ 
for given $\bf \in \ell_2(\checkbcJ)$.
{We anticipate that one does \emph{not} expect the usually dramatic effect of a squared condition number for $\bB^\top \bB$ since $\bB$ is wavelet-preconditioned, see below.}

\subsection{AWGMs for normal equations}

Now we {investigate} if the reformulation of \eqref{eq:equiv_ell2_prob} in terms of \eqref{eq:normal_equations} addresses the issues mentioned in Section \ref{subsec:extension_idealawgm_parabolic}.

\paragraph*{\emph{(1) Symmetry and positive definiteness}}
Obviously, $\bB^\top \bB$ is symmetric.
Moreover, by \eqref{eq;bounds_bB_bBinv}, it is also positive definite and it holds that
\begin{equation} \label{eq:norm_estim_bBtopbB}
   \| \bB^\top \bB \| \leq \| \bB \|^2, \quad    \| (\bB^\top \bB)^{-1} \| \leq \| \bB^{-1} \|^2,
\end{equation}
{hence} $\kappa(\bB^\top \bB)  \leq \| \bB \|^2  \| \bB^{-1} \|^2$.
Thus, we  consider $\bC \bu = \bg$ with $\bC = \bB^\top \bB$, $\bg = \bB^\top \bf$ and $\triplenorm \cdot \triplenorm^2 := \eval{\,\cdot\,}{\bB^\top \bB \, \cdot\,}$ and use Algorithm \ref{alg:exact_awgm}.

\paragraph*{\emph{(2) Bulk chasing and residual computation}}
Instead of considering the residual in $\ell_2(\checkbcJ)$, we now {obtain}  $\bB^\top (\bf - \bB \bu_{\hatbLambda_k}) \in \ell_2(\hatbcJ)$ with error estimator  $\rho_k:=\| \bB^\top (\bf - \bB \bu_{\hatbLambda_k})\|_{\ell_2(\hatbcJ)}$.
In analogy to \eqref{eq:estimate_elliptic_problem} and \eqref{eq:estimate_nonelliptic_problem}, we infer that
\begin{equation} \label{eq:estimate_leastsquares_problem}
  \| \bB \|^{-2} 
  \rho_k 
  \leq  \|   \bu -  \bu_{\hatbLambda_k} \|_{\ell_2(\hatbcJ)} 
  \leq \| \bB^{-1} \|^2 
  \rho_k.
\end{equation}
In this setting, the residual $\bf - \bB \bu_{\hatbLambda_k}$ from \eqref{eq:estimate_nonelliptic_problem} is also referred to as \emph{primal residual} whereas $\bB^\top(\bf - \bB \bu_{\hatbLambda_k})$ is called \emph{dual residual}.
Observe that this kind of residual allows for a bulk chasing strategy as used in line 5 of \textbf{IDEALIZED--AWGM}.

\paragraph*{\emph{(3) Well-posedness}}
With ${}_{\hatbLambda} \bB^\top$ and $\bB_{\hatbLambda}$ defined in \eqref{eq:bB_finite}, we get $(\bB^\top \bB)|_{\hatbLambda \times \hatbLambda} = \bbBB{\hatbLambda}{\hatbLambda}$ so that \eqref{eq:galerkin_system_elliptic} for general $\hatbLambda \subset \hatbcJ$ with $\bC = \bB^\top \bB$ now reads as follows:
\begin{equation} \label{eq:galerkin}
  \text{Find }\bu_{\hatbLambda} \in \ell_2(\hatbLambda): \qquad   \bbBB{\hatbLambda}{\hatbLambda} \bu_{\hatbLambda} = {}_{\hatbLambda}\bB^\top \bf, \qquad {}_{\hatbLambda}\bB^\top \bf \in \ell_2(\hatbLambda).
\end{equation}
Observe that the unique solution $\bu_{\hatbLambda}  = \argmin_{\bv_{\hatbLambda} \in \ell_2(\hatbLambda)} \| \bB_{\hatbLambda} \bv_{\hatbLambda} - \bf \|_{\ell_2(\checkbcJ)}^2$ to \eqref{eq:galerkin} can also be characterized as the solution of a least-squares problem. Moreover, the Galerkin problem \eqref{eq:galerkin} is  \emph{uniformly well-posed}.
Since $\bB^\top \bB$ is s.p.d., we infer from \eqref{eq:norm_estim_bBtopbB} that
$\| \bbBB{\hatbLambda}{\hatbLambda} \| \leq \| \bB \|^2$ as well as  $\|  ( \bbBB{\hatbLambda}{\hatbLambda} )^{-1} \| \leq \| \bB^{-1} \|^2$ for all $\hatbLambda \subseteq \hatbcJ$.
In particular, the condition number $\kappa(\bbBB{\hatbLambda}{\hatbLambda})$ is bounded \emph{independently} of $\hatbLambda$.


\begin{remark}
Obviously, neither the residual in \eqref{eq:estimate_leastsquares_problem} nor the solution $\bu_{\hatbLambda_k}$ of \eqref{eq:galerkin} can be computed exactly since the involved matrices are of infinite dimension. In order to obtain an implementable scheme, we work with an {approximation} $\bw_{\hatbLambda_k}$ to $\bu_{\hatbLambda_k}$ and an approximate residual $\res_k$ to $\bB^\top (\bf - \bB \bu_{\hatbLambda_k})$.
This will be discussed next. 
\end{remark}


\section{An implementable space-time adaptive wavelet Galerkin method} \label{sec:mt_awgm}
Now we describe the quasi-optimal ({in terms of} \ref{item:conv_rate} and \ref{item:comp_work})  AWGM for the numerical solution of \eqref{eq:normal_equations}  and call it \textbf{LS--AWGM} (\emph{least squares adaptive wavelet Galerkin method}), see Algorithm \ref{alg:ls_awgm}. 
We first describe the required subroutines.
We assume that $\bu\in\cA^s$ and denote by $\bw_{\hatbLambda}$ an approximate solution to \eqref{eq:galerkin}.

\begin{enumerate}[label=(\textbf{RES}), ref=(\textbf{RES}),leftmargin=35pt]
  \item \label{item:Res} \emph{Approximate residual:} For a given relative tolerance $0 < \omegals < 1$, the output $\resNE$ of $\textbf{RESIDUAL}[\bw_{\hatbLambda},\omegals]$ should satisfy
\begin{equation} \label{eq:req_resNE}
     \| \bB^\top(\bf - \bB \bw_{\hatbLambda})  - \resNE \|_{\ell_2(\hatbcJ)} \leq \omegals \cdot \nu, \quad \nu := \| \resNE \|_{\ell_2(\hatbcJ)},
\end{equation}
  and the associated computational cost is of order $\cO(\#\hatbLambda + \nu^{-1/s} \| \bu \|^{1/s}_{\cA^s})$.
\end{enumerate}
\begin{enumerate}[label=(\textbf{GAL}), ref=(\textbf{GAL}),leftmargin=37pt]
  \item \label{item:Gal} \emph{Approximate Galerkin problem:} For a given relative tolerance $0 < \gammals <1$, the output $\bw_{\hatbLambda}$ of $\textbf{GALSOLVE}[\hatbLambda, \overline{\bw}_{\hatbLambda}, \gammals \cdot \nu]$ should satisfy
  \begin{equation} \label{eq:req_Gal}
       \|  {}_{\hatbLambda}\bB^\top ( \bf -   \bB_{\hatbLambda} \bw_{\hatbLambda} ) \|_{\ell_2(\hatbLambda)} \leq \gammals \cdot \nu,
  \end{equation}
  where $\nu$ is defined in \eqref{eq:req_resNE} and the associated computational cost is of order $\cO(\#\hatbLambda + \nu^{-1/s}\| \bu \|^{1/s}_{\cA^s})$.
  Moreover, we assume that we are given an initial value $\overline{ \bw}_{\hatbLambda}$ satisfying $\|  {}_{\hatbLambda}\bB^\top ( \bf -   \bB_{\hatbLambda} \overline{\bw}_{\hatbLambda} ) \|_{\ell_2(\hatbLambda)} \leq (1+\gammals) \cdot \nu$.
\end{enumerate}
\begin{enumerate}[label=(\textbf{EXP}), ref=(\textbf{EXP}),leftmargin=37pt]
   \item \label{item:Exp} \emph{Approximate bulk chasing:} For a given parameter $0 < \delta <1$, the output $\underline{\hatbLambda} \subset \hatbcJ$ of  $\textbf{EXPAND}[\hatbLambda, \resNE,\delta]$ should satisfy
   \begin{equation} \label{eq:approx_bulk_chasing}
   \underline{\hatbLambda} \supset \hatbLambda, \qquad \| \bR_{\underline{\hatbLambda}} \resNE \|_{\ell_2(\underline{\hatbLambda})} \geq \delta \| \resNE \|_{\ell_2(\hatbcJ)},
   \end{equation}
   and, up to some absolute multiple, $\underline{\hatbLambda}$ is minimal among all sets that satisfy \eqref{eq:approx_bulk_chasing}.
   The computational cost of this routine is of order $\cO(\#\hatbLambda + \#\supp \resNE)$.
\end{enumerate}


\begin{algo} 
\caption{[$\bu_\eps$] = \textbf{LS--AWGM}[$\eps$, $\hatbLambda_1,\nu_0$]
	\label{alg:ls_awgm}}
\algsetup{indent=2em}
\begin{algorithmic}[1]
\REQUIRE Target tolerance $\eps$, finite index set $\hatbLambda_1 \subset \hatbcJ$ and tolerance $\nu_0 \eqsim \| \bB^\top \bf \|_{\ell_2(\hatbcJ)}$.
\ENSURE $\delta, \omegals, \gammals$ with \mbox{$\omegals \in (0, \delta)$}, \mbox{$\frac{\delta+\omegals}{1-\omegals} < \kappa(\bB^\top \bB)^{-\frac12}$},\newline $\gammals \in (0,\frac{(1-\omegals)(\delta-\omegals)}{1+\omegals} \kappa(\bB^\top \bB)^{-1})$.
\STATE Set $\bw_{\hatbLambda_0} := 0$.
\FOR{$k=1,2,\ldots$}
\STATE $\bw_{\hatbLambda_{k}} := \textbf{GALSOLVE}[\hatbLambda_{k},\bw_{\hatbLambda_{k-1}}, \gammals \cdot \nu_{k-1}]$.
\STATE $\resNE_k := \textbf{RESIDUAL}[\bw_{\hatbLambda_k},\omegals]$ and set $\nu_k := \| \resNE_k \|_{\ell_2}$.
\STATE \textbf{if} $\nu_k \leq \eps / \| \bB^{-1} \|^2$ \textbf{then} \textbf{return} $\bu_\eps := \bw_{\hatbLambda_k}$.
\STATE $\hatbLambda_{k+1} := \textbf{EXPAND}[\hatbLambda_k, \resNE_k, \delta]$
\ENDFOR
\end{algorithmic}
\end{algo}

In analogy to Proposition \ref{prop:conv_idealized_awgm}, we have the following result for \textbf{LS--AWGM} which is a direct consequence of \cite[Proposition 4.2 \& Theorem 4.1]{Stevenson:2009}:

\begin{theorem} [{\cite{Gantumur:2007,Stevenson:2009}}]
Let the assumptions on \emph{\ref{item:Res}}, \emph{\ref{item:Gal}} and \emph{\ref{item:Exp}} {and the requirements on $\delta, \omegals, \gammals$ from Algorithm \ref{alg:ls_awgm}} hold.
Then, the {iterates $\bw_{\hatbLambda_k}$} produced by  \textbf{\emph{LS--AWGM}} satisfy $\triplenorm \bu - \bw_{\hatbLambda_k} \triplenorm \leq \rho^{k/2} \triplenorm \bu \triplenorm$ where $\rho := 1- (\frac{\delta-\omegals}{1+\omegals})\kappa(\bB^\top\bB)^{-1} + \frac{\gammals^2}{(1-\omegals)^2}\kappa(\bB^\top \bB) < 1$  and {the output $\bu_\eps$ satisifes} $\| \bu - \bu_\eps \|_{\ell_2(\hatbcJ)} \leq \eps$.
If, moreover, $\bu \in \cA^s$ for some $s>0$, it holds for $\cN_k := \# \hatbLambda_k$ that
\begin{equation} \label{eq:quasi_optimal_lsawgm}
     \| \bu - \bw_{\hatbLambda_k} \|_{\ell_2(\hatbcJ)} \lesssim \| \bu \|^{1/s}_{\cA^s} \cN_k^{-s} , \quad \# \supp \bu_\eps \lesssim \eps^{-1/s} \| \bu \|^{1/s}_{\cA^s}.
\end{equation}
\end{theorem}


Note that a realization of \textbf{EXPAND} can easily {be} obtained by an \emph{approximate} sorting of $\resNE$ and a subsequent thresholding (e.g.\ \cite[p.\ 569]{Stevenson:2009}).
Possible realizations of the routines \textbf{RESIDUAL} and \textbf{GALSOLVE} that are based on so-called \textbf{APPLY}-routines (i.e., an adaptive, column-wise approximation of $\bB$ and $\bB^\top$) have been discussed in  \cite{Schwab:2009,Stevenson:2009}. 
We shall focus on a multitree approach which has been shown to outperform \textbf{APPLY}-based AWGMs in elliptic settings (see \cite{Kestler:2012d}).

\section{A multitree implementation} \label{sec:impl}


\subsection{Tree and multitree structured index sets} \label{subsec:trees_and_mtrees}
Let $\Psi = \{ \psi_\lambda : \lambda \in \cJ\}$ be a univariate uniformly local, piecewise polynomial wavelet basis as in \eqref{eq:Psi}.

\begin{definition} \label{def:tree}
A subset $\Lambda \subset \cJ$ is called a \emph{tree} if for any $\lambda \in \Lambda$ with $|\lambda|>0$ {it holds that
$\supp\,\psi_\lambda \subset \bigcup_{\mu\in\Lambda; |\mu|=\lambda-1} \supp\,\psi_\mu$.}
\end{definition}


It holds for all $\lambda,\mu \in \Lambda$ with
$|\mu| = |\lambda|-1$ and $|\supp \psi_\lambda \cap \supp \psi_\mu| > 0$ that $S_\mu \supset S_\lambda$, where
\begin{equation}
 S_\mu := \{ x \in \Omega: \dist(x,\supp \psi_\mu) \leq \rmD{\Psi}2^{-|\mu|} \}, \;\;  \rmD{\Psi} := \sup_{\lambda \in \cJ} 2^{|\lambda|} \diam (\supp \psi_\lambda). \hspace{-1ex}
\end{equation}

Let us now consider a tensor product wavelet basis $\bPsi = \{ \bpsi_{\lambda} : \blambda \in \bcJ \} \in \{ \hatbPsi, \checkbPsi\}$ with $\hatbPsi$ and $\checkbPsi$ as in \eqref{eq:hatPsiXL2} and \eqref{eq:hatPsiYL2}.
The extension of Definition \ref{def:tree} then reads:

\begin{definition}[{\cite{Kestler:2012c}}]\label{def:multitree}
An index set $\bLambda \in \bcJ$ is called a \emph{multitree} {if} for all $i \in \{ 0,\ldots,n\}$ {and all indices  $\mu_j \in \cJ^{(j)}$ for  $j \neq i$,  the index set
\begin{equation} \label{eq:tree_from_multitree}
   \Lambda^{(i)} := \{ \lambda_i \in \cJ^{(i)} :  (\mu_0, \ldots,\mu_{i-1},\lambda_i,\mu_{i+1},\ldots,\mu_n) \in \bLambda\}    \subset \cJ^{(i)}
\end{equation}
is either the empty set or a tree in the sense of Definition \ref{def:tree}}.
\end{definition}

Loosely speaking, a multitree $\bLambda \in \bcJ$ is ``when frozen in any $n$ coordinate directions, a tree in the remaining coordinate'' (see \cite[\S3.1]{Kestler:2012d}).

\begin{remark}
Note that quasi-optimality of \textbf{LS--AWGM} is maintained if $\hatbLambda_k$ {are required to be} multitrees (cf.\ \cite{Kestler:2012d}).
The only modification is to replace the \emph{unconstrained} nonlinear approximation space $\cA^s$ (see \eqref{eq:cAs}) by the \emph{constrained} approximation space
$\cAmtree^s := \{ \bv \in \ell_2(\hatbcJ) : \| \bv \|_{\cAmtree^s} < \infty \}$, 
where $\| \bv \|_{\cAmtree^s} := \sup_{\eps >0} \eps \cdot \big[  \min\{ \cN \in \N_0: \| \bv - \bv_\cN \|_{\ell_2(\hatbcJ)} \leq \eps \; \wedge \; \supp \bv_{\cN} \text{ is a multitree} \} \big]^s$.
This means we only allow {those} $\bv_\cN$ that are supported on a multitree.
\end{remark}

The reason for using trees and multitrees for solving linear operator equations instead of {arbitrary} index sets lies in the much more efficient evaluation of system matrices which we explain next.
{Moreover, tree and multitree-structured index sets are crucial {ingredients} for the evaluation of non-linearities in both tensor product settings (e.g.\ \cite{Stevenson:2011}) and non-tensor product settings  (e.g.\ \cite{Cohen:2003}). }

\subsection{Fast evaluation of {tensor product} system matrices}

We assume that for some $M \in \N$, there exist univariate bilinear forms $b^{(i)}_m$ such that
\begin{equation} \label{eq:prod_struc_b}
  \bB =  \bD^{\cY} \Big[   \sum_{m=1}^M \prod_{i=0}^n b_m^{(i)}(\hatPsi^{(i)},\checkPsi^{(i)}) \Big] \bD^{\cX}=  \bD^{\cY} \Big[ \sum_{m=1}^M \bigotimes_{i=0}^n \vecB^{(i)}_m \Big] \bD^{\cX},
\end{equation}
where $\vecB^{(i)}_m := b^{(i)}_m(\hatPsi^{(i)},\checkPsi^{(i)})$ for $i=0,\ldots,n$ and $m=1,\ldots,M$. {This means that $\bB$ is a preconditioned sum of tensor product bilinear forms. As we shall see below, this form holds true for a large class of operators.} 
{Moreover}, we shall {always} assume that $b^{(i)}_m$ are \emph{local} in the sense that $b^{(i)}_m(w,v) = 0$ whenever $|\supp v \cap \supp w| = 0$. 
%
%
The special structure of  $\bB$  can  be used to {efficiently} realize the application of  $\bbB{\checkbLambda}{\hatbLambda}$ to a vector $\bv_{\hatbLambda} \in \ell_2(\hatbLambda)$ for finite multitrees $\hatbLambda \in \hatbcJ$ and $\checkbLambda \in \checkbcJ$.
As described in \cite{Kestler:2012c}, this can be realized in linear complexity, i.e., $\cO(\# \hatbLambda + \#\checkbLambda)$ by using a \emph{separation} of $\bbB{\checkbLambda}{\hatbLambda}$ into \emph{unidirectional} operations and an efficient \emph{tree-based} application of  unidirectional operations.
These principles are also known from sparse grid algorithms (see, e.g., \cite{Zenger:1991,Bungartz:2004}).

We recall the Kronecker product of two general {(possible  bi-infinite)} matrices $\vecA^{(1)}, \vecA^{(2)}$ and identity matrices $\vecId{}^{(1)}, \vecId{}^{(2)}$ of appropriate dimension:
\begin{equation} \label{eq:kronecker_splitting}
\kern-5pt
\vecA^{(1)} \otimes \vecA^{(2)} 
=\kern-3pt \big[\vecA^{(1)} \otimes \vecId{}^{(2)}\big] \circ \big[\Id{}^{(1)} \otimes \vecA^{(2)} \big] 
=\kern-3pt \big[\vecId{}^{(1)} \otimes \vecA^{(2)} \big] \circ \big[\vecA^{(1)} \otimes \vecId{}^{(2)}\big].
\end{equation}
Then we split $\vecB^{(i)}_m = \vecL^{(i)}_m + \vecU^{(i)}_m$ into a (stricly) lower $\vecL^{(i)}_m := \big[ (\vecB^{(i)}_m)_{\lambda,\mu}]_{|\lambda|>|\mu|}$ and an upper triangular matrix  $\vecU^{(i)}_m := \big[ (\vecB^{(i)}_m)_{\lambda,\mu}]_{|\lambda|\leq|\mu|}$.
With \eqref{eq:kronecker_splitting}, it can {then} be shown that there exist multitrees $\overline{\bXi}$ and $\underline{\bXi}$ such that we have the following \emph{equivalent} representation of  $\bbB{\checkbLambda}{\hatbLambda}$,
\begin{align*}
 \bD^{\cY} \bigg[
 	  \sum_{m=1}^M \underbrace{ \bR_{\checkbLambda} \big[\vecId{}^{(0)} \otimes \vecB_m^{(1)} \otimes \cdots \otimes \vecB_m^{(n)} \big] \bE_{\overline{\bPi}} }_{=: \textrm{(I)}} \circ  \underbrace{ \bR_{\overline{\bPi}}\big[ \vecU_m^{(0)} \otimes \vecId{}^{(1)} \otimes \cdots \otimes \vecId{}^{(n)}  \big] \bE_{\hatbLambda} }_{=: \textrm{(II)}} \\
	\kern-11pt
    +    \sum_{m=1}^M \underbrace{ \bR_{\checkbLambda}  \big[\vecL_m^{(0)} \otimes \vecId{}^{(1)} \otimes \cdots \otimes \vecId{}^{(n)} \big]  \bE_{\underline{\bPi}} }_{=: \textrm{(III)}} \circ  \underbrace{ \bR_{\underline{\bPi}} \big[ \vecId{}^{(0)} \otimes \vecB_m^{(1)} \otimes \cdots \otimes \vecB_m^{(n)} \big] \bE_{\hatbLambda}  }_{=: \textrm{(IV)}} \bigg] \bD^{\cX}.
\end{align*}
It holds $\#\overline{\bXi} + \#\underline{\bXi} \lesssim \# \checkbLambda + \# \hatbLambda$. 
The application of (II), (III) (and (I), (IV) for $n=1$) is referred to as \emph{unidirectional operation} as  only the application of the univariate matrices
$\vecL^{(0)}_m|_{\checkLambda^{(0)} \times \hatLambda^{(0)}}$, $\vecU^{(0)}_m|_{\checkLambda^{(0)} \times \hatLambda^{(0)}}$ and $\vecB^{(1)}_m|_{\checkLambda^{(1)} \times \hatLambda^{(1)}}$ ($n=1$) is required. 
%
Due to the tree structure, these tasks can be realized in linear complexity despite the fact that neither of the matrices $\vecL^{(0)}_m$, $\vecU^{(0)}_m$ or $\vecB_m^{(1)}$ is sparse {in general} (see \cite[\S2]{Kestler:2012c}).
For $n>2$, the remaining {parts (I) and (IV)} can be treated recursively by applying the same procedure to  $\vecB_m^{(1)} \otimes \cdots \otimes \vecB_m^{(n)}$. 

%
%

\begin{theorem}[{\cite[Theorem 3.1]{Kestler:2012c}}] \label{thm:mv_mt}
Let $\cA$ be a linear differential operator with polynomial coefficients and {let} $\hatbLambda \subset \hatbcJ$, $\checkbLambda \in \checkbcJ$ be multitrees.
Then, for any $\bv_{\hatbLambda} \in \ell_2(\hatbLambda)$, the product $\bbB{\checkbLambda}{\hatbLambda} \bv_{\hatbLambda}$ can be {computed} in $\cO(\# \hatbLambda + \# \checkbLambda)$ operations.
\end{theorem}

\begin{remark}
{If $\cA$ is a linear differential operator with polynomial coefficients,  $\bB$ has the form  \eqref{eq:prod_struc_b}.} 
Furthermore, all matrices  can be applied in linear complexity {if} $\hatLambda^{(i)}$ and $\checkLambda^{(i)}$ are trees (cf.\ \cite[\S2]{Kestler:2012c}).
\end{remark}


\subsection{\textbf{RESIDUAL:} Multitree residual approximation}
We need to  approximate the residual $\bB^\top (\bf - \bB \bw_{\hatbLambda})$ by a residual of type
$\bbBtr{\hatbXi}{\checkbXi} \;\big( \bf_{\checkbXi} - \bbB{\checkbXi}{\hatbLambda} \bw_{\hatbLambda} \big)$.

\subsubsection{Primal residual}
We first {recall} the approximation of the primal residual.

\begin{theorem}[{\cite{Kestler:2012d}}] \label{thm:primal_res}
Let $0< \omega <1$, {let} $\cA$ be a differential operator with polynomial coefficients and {let} $\bu \in \cAmtree^s$ for some $s>0$.
Then,  for all finite multitrees $\hatbLambda \subset \hatbcJ$ and {all} $\bw_{\hatbLambda} \in \ell_2(\hatbLambda)$, there exists a multitree $\checkbXi = \checkbXi(\hatbLambda,\omega) \subset \checkbcJ$ such that $\# \checkbXi \lesssim \# \hatbLambda + \nu^{-1/s}$ with $\nu := \| \res \|_{\ell_2(\checkbcJ)}$,   $\bf_{\checkbXi} := \bR_{\checkbXi} \bf$ and
\begin{equation} \label{eq:error_estim_primal_res}
    \| (\bf - \bB \bw_{\hatbLambda}) - \res \|_{\ell_2(\checkbcJ)} \leq \omega \| \res \|_{\ell_2(\checkbcJ)}, \quad \res := \bf_{\checkbXi} - \bbB{\checkbXi}{\hatbLambda} \bw_{\hatbLambda}.
\end{equation}
\end{theorem}

\begin{remark}
{Due} to the multitree structure of $\hatbLambda$ and $\checkbXi$, the computational cost for \emph{computing} $\res$ {is} $\cO(\# \hatbLambda + \nu^{-1/s})$ if an entry $\bf_{\blambda}$ of $\bf = (\bf_{\blambda})_{\blambda \in \checkbcJ}$ can be computed \emph{exactly at unit cost}, which is e.g.\ the case if $f$ is a {(piecewise)} polynomial. 
If this assumption is not met, replace $\bf$ by some $\bf_\eps$ with $\| \bf - \bf_\eps \|_{\ell_2(\checkbcJ)} \leq \varepsilon$ and $\# \supp \bf_\eps \lesssim \eps^{-1/s}$ which is possible if $f$ is sufficiently (piecewise) smooth (see \cite[\S3.4]{Kestler:2012d}).
\end{remark}

\subsubsection{Dual residual}\label{Sec:RESDualRes}
We may now follow \cite[\S1.1]{Kestler:2012d} using a wavelet compression of  $\bB$ and $\bB^\top$. 
If $\cA$ is {a} linear differential operator with polynomial coefficients, it can be shown that for any $0< \eta <1$, there exists $\bB_\eta: \ell_2(\hatbcJ) \to \ell_2(\checkbcJ)$ such that
\begin{equation} \label{eq:compr}
   \| \bB - \bB_\eta \| \leq \eta, \quad \| \bB^\top - \bB^\top_\eta \|  \leq \eta,
\end{equation}
where the number of nonzeros in each row and each column of $\bB_\eta$ are of order $\cO(\eta^{-1/\sstar})$ for some $\sstar > \smax$, \eqref{eq:smax}. {This means that} $\bB$ is $\sstar$-admissible  (see \cite{Schwab:2009}). 
%
Assuming   that $\eta$ is chosen sufficiently small so that $\bB_\eta$ and $\bB^\top_\eta$ are boundedly invertible, we obtain the estimate (see Proposition \ref{prop:dual_res})
\begin{equation} \label{eq:error_estim_dual_res}
   \| \bB^\top (\bf - \bB \bw_{\hatbLambda}) - \bB_\eta^\top \res \|_{\ell_2(\hatbcJ)} \leq \omegals \| \bB^\top_\eta \res \|_{\ell_2(\hatbcJ)},
\end{equation}
for $\omegals  = (\eta{\textstyle \frac{1}{1-\omega} } + (\| \bB \| + \eta) \omega) \| \bB^{-1}_\eta \|$ so that $\omegals \to 0$ as $\omega\to 0$ and $\eta \to 0$. 
Even though $\bB_\eta$ and $\bB_\eta^\top$ are sparse (for fixed $\eta$), the application of these matrices to finite vectors can be computationally expensive since the product structure of $\bB$ in \eqref{eq:prod_struc_b} cannot {be} exploited.
Unfortunately, the approximate residual $\bB^\top_\eta \res$ is \emph{not} necessarily supported on a {multitree. Hence,} we define the multitree-based residual
\begin{equation} \label{eq:resNE}
   \resNE := \bbBtr{\hatbXi}{\checkbXi} \;\big( \bf_{\checkbXi} - \bbB{\checkbXi}{\hatbLambda} \bw_{\hatbLambda} \big) =  \bbBtr{\hatbXi}{\checkbXi} \; \res
\end{equation}
{such that} $\| \bB^\top (\bf - \bB \bw_{\hatbLambda}) - \resNE \|_{\ell_2(\hatbcJ)} \leq \omegals \| \resNE \|_{\ell_2(\hatbcJ)}$ where $\hatbXi$ is the smallest multitree containing $\supp \bB^\top_\eta \res$. The residual computation requires $\cO(\#\hatbXi + \#\checkbXi)$ operations.

\begin{remark} \label{rem:choiceofindexsets}
Theorem \ref{thm:primal_res} only ensures the existence of an appropriate multitree $\checkbXi$ but does not give any information on its explicit construction.
The same holds true for $\hatbXi$.
In Section \ref{Sec:TestSetChoice}, we will discuss how we can construct the multitrees $\checkbXi$ and $\hatbXi$  \emph{without} setting up the compressed matrix $\bB^\top_\eta$ so that $\resNE$ from \eqref{eq:resNE} satisfies \ref{item:Res}.
Furthermore, numerical experiments in Section \ref{sec:numerics} indicate appropriate choices of $\checkbXi$ and $\hatbXi$ with \emph{preferably small} cardinalities and \emph{optimal balancing} of the error arising from the approximations of the primal (see \eqref{eq:error_estim_primal_res}) and dual residual (see \eqref{eq:error_estim_dual_res}).
\end{remark}

\subsection{\textbf{GALSOLVE:} Multitree solution of finite-dimensional least squares problems}
Concerning the numerical solution of the least squares problem \eqref{eq:normal_equations}, the approach proposed in \cite{Stevenson:2009, Schwab:2009} consists of replacing $\bbBB{\hatbLambda}{\hatbLambda}$  by a \emph{sparse} approximation ${}_{\hatbLambda}[\bB_\eta^\top \bB_\eta]_{\hatbLambda} := \bR_{\hatbLambda}[\bB_\eta^\top \bB_\eta] \bE_{\hatbLambda}$ satisfying $\| \bbBB{\hatbLambda}{\hatbLambda} - {}_{\hatbLambda}[\bB_\eta^\top \bB_\eta]_{\hatbLambda}\| \lesssim \eta$.
In analogy to \eqref{eq:normal_equations}, we consider:
\begin{equation} \label{eq:gal_eta_version}
    \text{Find } \bu_{\eta,\hatbLambda} \in \ell_2(\hatbLambda):  \qquad  {}_{\hatbLambda}[\bB_\eta^\top \bB_\eta]_{\hatbLambda} \bu_{\eta,\hatbLambda}  =  \bR_{\hatbLambda} \bB_{\eta}^\top \bf_{\checkbLambda}.
\end{equation}
Indeed, under the assumption that $\eta$ is sufficiently small, $\kappa({}_{\hatbLambda}[\bB_\eta^\top \bB_\eta]_{\hatbLambda})$ is bounded independently of $\hatbLambda$ (see Appendix \ref{sec:appendix}).
In particular, there exist algorithms based on linear iterative solvers like {the conjugate gradient (cg) method} that {approximate}  \eqref{eq:gal_eta_version} such that $\|\bu_{\hatbLambda} - \bu_{\eta,\hatbLambda} \|_{\ell_2(\hatbLambda)} \| \lesssim \eta$ and \ref{item:Gal} is satisfied. 
Similar to the residual approximation, the disadvantage of this approach is that we cannot use the fast matrix-vector multiplication w.r.t.\ multitrees. 
To this end, we intend to compute $\bw_{\hatbLambda}$ as an approximate solution of the problem: 
\begin{equation} \label{eq:normalequations_heuristic}
 \text{Find } \bx_{\hatbLambda} \in \ell_2(\hatbLambda):  \qquad  \bbBtr{\hatbLambda}{\checkbLambda} \bbB{\checkbLambda}{\hatbLambda} \bx_{\hatbLambda}  =  \bbBtr{\hatbLambda}{\checkbLambda} \bf_{\checkbLambda}.
\end{equation}
We could choose $\checkbLambda$ as the smallest multitree that contains $\supp \bB_\eta \bv_{\hatbLambda}$ for all $\bv_{\hatbLambda} \in \ell_2(\hatbLambda)$.
However, this is \emph{not} an implementable approach. 
Hence, we are concerned with the question how the multitree $\checkbLambda$ can be constructed in dependency of $\hatbLambda$ such that (1) the condition number of $\bbBtr{\hatbLambda}{\checkbLambda}$ is uniformly bounded and (2) an approximate solution $\bw_{\hatbLambda}$ to \eqref{eq:normalequations_heuristic} satisfies \ref{item:Gal}. This will be discussed in {Sections \ref{Sec:TestSetChoice} and \ref{sec:numerics}}. For \emph{fixed} multitrees, the solution of \eqref{eq:normalequations_heuristic}  can be computed e.g.\ with cg.

\subsection{Choice of index sets}\label{Sec:TestSetChoice}

The expansion $\hatbLambda_{k} \to \hatbLambda_{k+1}$ of the \emph{trial}  sets in Algorithm \ref{alg:ls_awgm} is based upon the residual  $\resNE_{k}$, but it is so not clear how to construct appropriate \emph{test} sets $\checkbLambda_{k}= \checkbLambda_{k}(\hatbLambda_{k})$. Similarly for the auxiliary  sets $\hatbXi_{k}$ and $\checkbXi_{k}$ required for \eqref{eq:resNE}: {While} the construction of the test sets $\checkbXi_{k}$ for the primal residual in a Galerkin setting has been investigated in \cite{Kestler:2012d}, there are {so far} no  results for good choices of $\checkbXi_{k}$ and $\hatbXi_{k}$ within a Petrov-Galerkin framework.

\subsubsection*{Choice of test sets $\checkbLambda_{k}$}
For a given index set $\hatbLambda_{k} \in \hatbcJ$, we have to ensure that the finite-dimensional test set $\checkbLambda_{k} \in \checkbcJ$ is large enough to ensure well-posedness. At the same time, for efficiency we would like to choose $\checkbLambda_{k} \in \checkbcJ$ as small as possible. We describe a corresponding iteration. 
As initial  sets $\hatbLambda_{0}$, $\checkbLambda_{0}$, we follow \cite[\S6.2]{Andreev:2012} 
\begin{align}
\hatbLambda_{0} 
	= \hatbLambda_{SG, J}
	&:= \{\blambda \in \hatbcJ: |\blambda| \leq J \}, \label{eq:SGhatLambda}\\
\checkbLambda_{0} 
	= \checkbLambda_{SG, J}
	&:= \{\blambda \in \checkbcJ: |\blambda| \leq J \text{ or } |\lambda_{0}| = J+1, 
		|\lambda_{i}| = 0, 1\le i\le n\}, \label{eq:SGcheckLambda}
\end{align}
where $|\blambda| := \sum_{i=0}^{n} |\lambda_{i}|$.\footnote{{We will also use $\hatbLambda_{SG, J}$ and $\checkbLambda_{SG, J}$ within a uniform sparse grid (SG) discretization.}} 
Such bases are provably stable, however, this only holds true for \emph{uniform} (full or sparse) discretizations. 
In later iterations, i.e. for adaptively constructed trial sets $\hatbLambda_{k}$, $k>0$, we propose the following (heuristic) choices:
\begin{enumerate}\renewcommand{\labelenumi}{(\roman{enumi})}
\item
 $\checkbLambda_{\textrm{Full}} = \textbf{FullStableExpansion}(\hatbLambda,\ell)$
 {is defined as}
\begin{align}
\quad\checkbLambda_{\textrm{Full}} &{:=} \bigl\{\blambda \in \checkbcJ: \exists\, \bmu \in \hatbLambda \text{ s.t. for all } j=0,\dots,n:  |\lambda_{j}| \leq |\mu_{j}| + \ell \label{Eq:FullStExp}\\
&\qquad
\text{ and } \dist\bigl(\supp\checkpsi_{\lambda_{j}}^{(j)},\supp\hatpsi_{\mu_{j}}^{(j)}\bigr) \leq \mathrm{D}_{\checkPsi^{(j)}}  2^{-|\lambda_{j}|}   \bigr\}. \notag
\end{align}
\item $\checkbLambda_{\textrm{Red}} = \textbf{ReducedStableExpansion}(\hatbLambda,\ell)$ {is a subset of $\checkbLambda_{\textrm{Full}}$ defined as}
\begin{align}
\checkbLambda_{\textrm{Red}} &:= \bigcup_{i=0}^{n} \bigl\{\blambda \in \checkbcJ: \exists\, \bmu \in \hatbLambda \text{ s.t. for all } j=0,\dots,n:  |\lambda_{j}| \leq |\mu_{j}| + \delta_{i,j}\,\ell  \label{Eq:RedStExp}\\[-3mm]
&\qquad\quad
\text{ and } \dist\bigl(\supp\checkpsi_{\lambda_{j}}^{(j)},\supp\hatpsi_{\mu_{j}}^{(j)}\bigr) \leq \mathrm{D}_{\checkPsi^{(j)}}  2^{-|\lambda_{j}|} \bigr\}. \notag
\end{align}\vspace{-2ex}
\item $\checkbLambda_{\textrm{Temp}} = \textbf{TemporalStableExpansion}(\hatbLambda,\ell)$: {consists of}  only \emph{temporal} higher level extensions, i.e., 
\begin{align}
 \quad\checkbLambda_{\textrm{Temp}} &:= \bigl\{\blambda \in \checkbcJ: \exists\, \bmu \in \hatbLambda \text{ s.t. for all } j=0,\dots,n:  |\lambda_{j}| \leq |\mu_{j}| + \delta_{0,j}\,\ell  \label{Eq:TempStExp}\\
&\qquad
\text{ and } \dist\bigl(\supp\checkpsi_{\lambda_{j}}^{(j)},\supp\hatpsi_{\mu_{j}}^{(j)}\bigr) \leq \mathrm{D}_{\checkPsi^{(j)}}  2^{-|\lambda_{j}|}\bigr\}. \notag
\end{align}
\end{enumerate}
{We refer} to \cite[Prop. 2]{Kestler:2012d} for a proof that the above index sets are indeed multitrees. An algorithmic realization is shown in Algorithm \ref{alg:stableexp}. 

\begin{algo} 
\caption{[$\checkbLambda$] = \textbf{FullStableExpansion}[$\hatbLambda$, $\ell$]
	\label{alg:stableexp}}
\algsetup{indent=2em}
\begin{algorithmic}[1]
\REQUIRE Finite index set $\hatbLambda \subset \hatbcJ$, expansion level $\ell \in \N$.
\STATE $\checkbLambda := \emptyset \subset \checkbcJ$.
\FOR{$\blambda =(\lambda_{0},\dots,\lambda_{n}) \in \hatbLambda$}
	\STATE Find all ``neighbours'' $\bmu  = (\mu_{0},\dots,\mu_{n}) \in \checkbcJ$ on the same level: 
	 \mbox{$\checkbLambda \gets \checkbLambda \cup \{\bmu \in \checkbcJ: |\mu_{i}| = |\lambda_{i}| , \ \supp \checkpsi_{\mu_{i}} \cap \supp \hatpsi_{\lambda_{i}} \neq 0  \ \forall\, i=0,\dots,n\}$}.
	\STATE Find all ``neighbours'' $\tilde\bmu  = (\tilde\mu_{1},\dots,\tilde\mu_{n}) \in \checkbcJ$ on the $\ell$ higher levels:\\
	\mbox{$\checkbLambda \gets \checkbLambda \cup \{\tilde\bmu \in \checkbcJ: |\tilde\mu_{i}| = |\lambda_{i}|+j, 1\leq j \leq \ell, \ \supp \checkpsi_{\tilde\mu_{i}} \cap \supp \hatpsi_{\lambda_{i}} \neq 0$} \mbox{$\hspace{1.8cm}\forall\, i=0,\dots,n\}$}.
	\STATE Complete $\checkbLambda$ to form a multitree in the sense of Definition \ref{def:multitree}.
\ENDFOR
\end{algorithmic}
\end{algo}

\subsubsection*{Choice of sets $\hatbXi_{k}$, $\checkbXi_{k}$}
The proposed index set reads
\begin{align}
\hatbXi_{k} &\ = \textbf{ReducedMultiTreeCone}(\hatbLambda,\ell) \notag\\
		&:= \bigcup_{i=0}^n \bigl\{\blambda \in \hatbcJ: \exists\, \bmu \in \hatbLambda_{k} \text{ s.t. for all } j=0,\dots,n:  |\lambda_{j}| \leq |\mu_{j}| + \delta_{i,j}\,\ell \label{eq:RedMTCone}\\[-3mm]
		&\qquad\qquad \text{ and } 
\dist\bigl(\supp\hatpsi_{\lambda_{j}}^{(j)},\supp\hatpsi_{\mu_{j}}^{(j)}\bigr) \leq \mathrm{D}_{\hatPsi^{(j)}} 2^{-|\lambda_{j}|} \bigr\} \notag
\end{align}
It was shown {in \cite{Kestler:2012c,Kestler:2012d}} that this index set for $\ell=1$ and the analogously defined {$\textbf{FullMultiTreeCone}(\hatbLambda,1)$} are adequate choices {for} an accurate approximation of the primal residual in the Galerkin setting, where $\hatbPsi^\cX=\checkbPsi^\cY$, $\hatbXi=\checkbXi$.

\begin{figure}[!ht]
\subfigure[\textbf{FullResConstruction}]{
\begin{tikzpicture} [node distance=1.7cm]
	  \tikzstyle{every node}=[font=\large]
		\node 										 (XLambda1) {$\hatbLambda_{k}$};
		\node [right=1mm of XLambda1] (XSubSet1) {\scalebox{1.5}{$\subset$}};
		\node [left=1mm of XLambda1] (XSpace){$\cX$:};
		\node [right=1mm of XSubSet1] (XXi1) {$ \hatbXi_{k}^{\textrm{tmp}}$};
		\node [above=-0.1cm of XSubSet1, align=center, font=\footnotesize, text width = 1.4cm] {\eqref{eq:RedMTCone}};
		\node [below of=XXi1, xshift=-0.2cm] (YXi1) {$\checkbXi_{k}$};
		\node [right of=XXi1] (XXi2) {$ \hatbXi_{k}^{\ }$};
		\node [below of=XSpace] (YSpace){$\cY$:};

		\draw [->, semithick, decorate, decoration={snake, segment length=7pt,amplitude=1pt, post length=2pt}] ($(XXi1.south)-(0.3,0)$) -- ($(YXi1.north)-(0.1,0)$)
			node[midway, left, xshift=-0.5mm, yshift=-1mm, text width=1.7cm, font=\footnotesize, inner sep=0pt]{\eqref{Eq:FullStExp}-\eqref{Eq:TempStExp}};
			
		\draw [->, semithick, decorate, decoration={snake, segment length=7pt,amplitude=1pt, post length=2pt}] ($(YXi1.east)-(0.1,-0.05)$) -- ($(XXi2.south)-(0.2,-0.1)$)
			node[midway, right, xshift=2.1mm, yshift=0.3mm, font=\footnotesize, inner sep=0pt]{\eqref{Eq:FullStExp}};
		\begin{pgfonlayer}{background}
			\node [fill=gray!30,fit=(current bounding box.north west) (current bounding box.south east)] {}; 
		\end{pgfonlayer}
	\end{tikzpicture}\label{Fig:XiConstr-Double}} \hspace{2mm}
\subfigure[\textbf{OptimResConstruction}]{
\begin{tikzpicture} [node distance=1.7cm]
	  \tikzstyle{every node}=[font=\large]
		\node 										 (XLambda1) {$\hatbLambda_{k}$};
		\node [left=1mm of XLambda1] (XSpace){$\cX$:};
		\node [right=1mm of XLambda1] (XSubSet1) {\scalebox{1.5}{$\subset$}};
		\node [right=1mm of XSubSet1] (XNabla1) {$ \hatbXi_{k}^{\textrm{tmp}}=: \hatbXi_{k}$};
		\node [above=-0.1cm of XSubSet1, align=center, font=\footnotesize, text width = 1.4cm] {\eqref{eq:RedMTCone}};
		\node [below of=XXi1, xshift=-0.2cm] (YXi1) {$\checkbXi_{k}$};
		\node [below of=XSpace] (YSpace){$\cY$:};
		\draw [->, semithick, decorate, decoration={snake, segment length=7pt,amplitude=1pt, post length=2pt}] ($(XXi1.south)-(0.3,0)$) -- ($(YXi1.north)-(0.1,0)$)
			node[midway, left, xshift=-0.5mm, yshift=-1mm, text width=1.7cm, font=\footnotesize, inner sep=0pt]{\eqref{Eq:FullStExp}-\eqref{Eq:TempStExp}};
		\begin{pgfonlayer}{background}
			\node [fill=gray!30,fit=(current bounding box.north west) (current bounding box.south east)] {}; 
		\end{pgfonlayer}
	\end{tikzpicture}\label{Fig:XiConstr-Simple}}
\caption{Constructions of index sets $\hatbXi_{k}$, $\checkbXi_{k}$ for residual approximation.}
\label{Fig:XiConstr}
\end{figure}
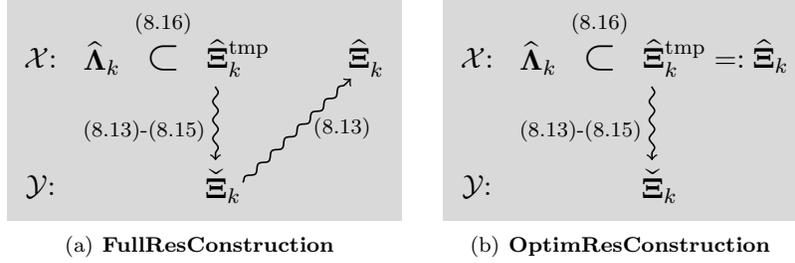

In our Petrov-Galerkin setting, we combine the multitree cone extension with the expansions \eqref{Eq:FullStExp}-\eqref{Eq:TempStExp}. More precisely, we {consider} the two variants \textbf{FullResConstruction} and \textbf{OptimResConstruction}, see Figure \ref{Fig:XiConstr}. 
{For the primal residual (i.e., in $\cY$), we expand $\hatbLambda_{k}$ to $\hatbXi_{k}^{\textrm{tmp}} = \textbf{ReducedMultiTreeCone}\bigl(\hatbLambda_{k},\ell\bigr)$ and obtain the desired $\checkbXi_{k}$ by one of the expansion variants in \eqref{Eq:FullStExp}-\eqref{Eq:TempStExp}.} 
%
%
For the dual residual {(in $\cX$), we consider two approaches. In the first one, shown in Fig.\ \ref{Fig:XiConstr-Double}, we take the  set $\checkbXi_{k}$ as above and set $\hatbXi_{k} = \textbf{FullStableExpansion}\bigl(\checkbXi_{k},\ell\bigr)$ (with obvious inverted roles of primal and dual basis).  Then, $\hatbXi_{k}$ is the smallest multitree containing $\supp \bB^\top_\eta \res_{k}$ for sufficiently small $\eta$. } 
{The second approach uses the by far smaller set  $\hatbXi_{k} = \hatbXi_{k}^{\textrm{tmp}}$ as indicated in Fig.\ \ref{Fig:XiConstr-Simple}, \cite{Kestler:2012d}.}

\section{{Numerical Experiments}} \label{sec:numerics}
We report numerical examples for  time-periodic problems of type \eqref{eq:problem_strongform}.
We focus on {the} \emph{stability of {the} arising normal equations}  \eqref{eq:normalequations_heuristic} in view of different choices for $\checkbLambda_k$. Moreover, we {numerically investigate} the \emph{quantitative behavior of approximate primal and dual residuals} in view of Remark \ref{rem:choiceofindexsets}. 
It is sufficient to consider the case $n=1$ (so that $\bOmega = (0,1)$), since we employ an $L_{2}(0,1)$-orthonormal \mbox{(multi-)}wavelet basis $\bSigma = \Sigma$ (see \eqref{eq:bSigma}) as in \cite{Rupp:Diss}, with $d_{x}=2$ and  homogeneous boundary conditions.
In this case, the Riesz constants  in \eqref{eq:rmc_V}, \eqref{eq:rmC_V} are \emph{independent} of $n$.
In particular, the condition numbers of $\bB^\top \bB$ and of $\bbBtr{\hatbLambda}{\checkbLambda} \bbB{\checkbLambda}{\hatbLambda}$ do not depend on $n$ so that the 1D case gives all relevant information. 
In \cite{Kestler:2012d}, it was shown numerically that the \emph{asymptotic} behavior of the multitree-based residual only differs by a constant depending on $n$ {from the unconstrained case}. 

We choose $\Theta^{\per}$ (see \eqref{eq:Thetaper}) as a collection of bi-orthogonal B-spline wavelets of order $d_{t}=\tilde d_{t}=2$ on the real line, periodized onto $[0,T]$, \cite{Urban:WaveletBook}. 
For $\Theta$ (see \eqref{eq:Theta}), we choose  bi-orthogonal B-spline wavelets  from \cite{Dijkema:Diss} with $d_{t}=\tilde d_{t}=2$.
{As further} parameters for the \textbf{LS-AWGM} {we choose} $\delta=0.7$\footnote{We have chosen a larger value for $\delta$ than required by Algorithm \ref{alg:ls_awgm} for efficiency reasons.}, $\gammals=0.01$ and, if not indicated differently, $\ell=1$ for the stable extensions from Section \ref{Sec:TestSetChoice}. 
%
%
%
%
%
We obtain {qualitatively} similar results for choosing  $\Sigma$ as in  \cite{Dijkema:Diss} 
for $d_{x}=\tilde d_{x} = 2$ even though they do not satisfy our assumptions.\footnote{Note that these bases cannot be normalized to be a Riesz basis of $H^{-1}(\Omega)$.} 

We also compare the \textbf{LS--AWGM} to a (uniform) sparse grid approach {(SG)}, i.e., to computing the solutions on a sequence of uniform finite-dimensional sets $\hatbLambda_{SG, J}$, $\checkbLambda_{SG, J}$, $J=0,1,\dots$, as in \eqref{eq:SGhatLambda}, \eqref{eq:SGcheckLambda}, {e.g.\ \cite{Zenger:1991,Bungartz:2004}}.

\subsection{Heat Equation} We consider the 1D-inhomogeneous heat equation \\
\begin{minipage}{0.68\textwidth}
\begin{equation*}
	\left\{
	\begin{aligned}
			u_{t} - u_{xx} &= f(t,x) \qquad \text{ on } \Omega = (0,1),\\
			u(t,0) &= u(t,1) \qquad \text{ for all } t \in [0,T],\\
			u(0,x) &= u(T,x)	 \qquad \!\text{on } \overline\Omega,
	\end{aligned}\right.
\end{equation*}
\end{minipage}\hfill
\begin{minipage}{0.3\textwidth}
%
%
%
\ \\[1ex]
\includegraphics[width=0.7\textwidth]{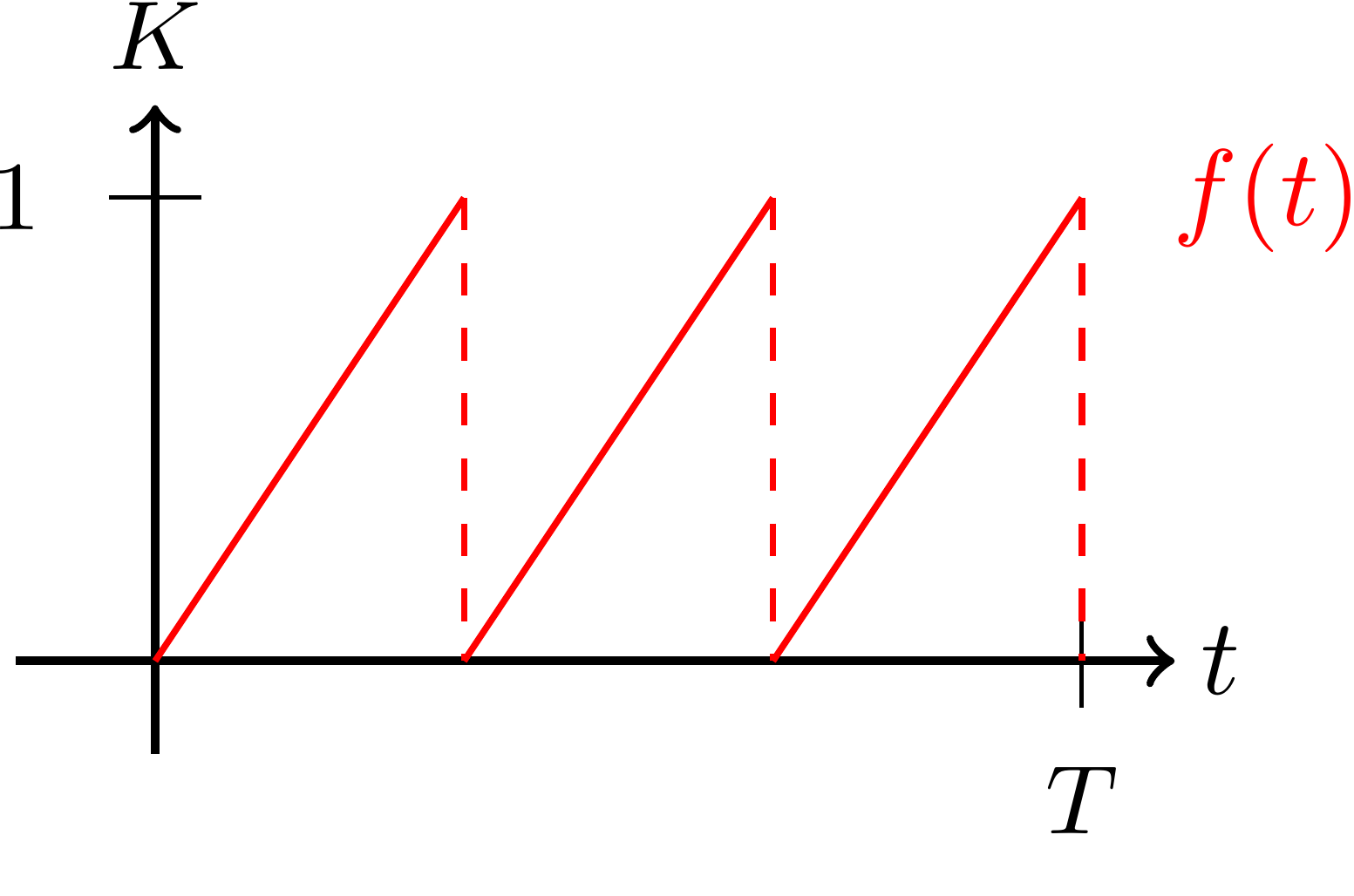}
\end{minipage}
with a {discontinuous} source function $f(t,x) \equiv f(t) := K\left(\tfrac{Nt}{T} - \lfloor \tfrac{Nt}{T} \rfloor \right)$, $N \in \N$, $K \in \R_{+}$. {Our figures  correspond to the choice $N=3$, $K=1$.}

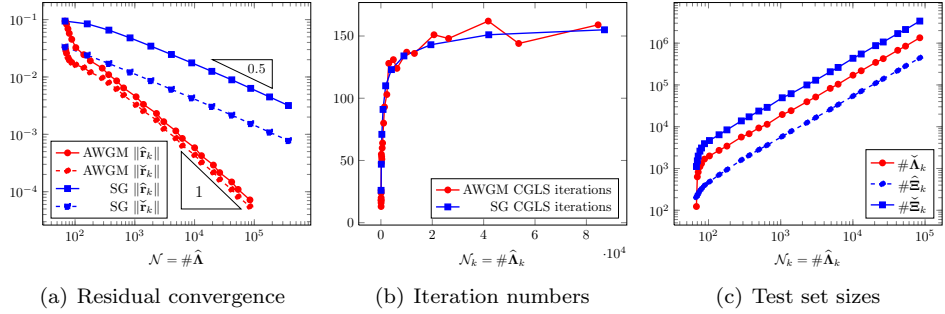
\begin{figure}[!ht]
	\subfigure[Residual convergence]{\begin{tikzpicture}[scale=0.52]
	  \begin{axis}
	  		[	xlabel={$\mathcal{N} = \#\hatbLambda$},
			xmode = log,
			ymode = log,
			ylabel near ticks,
			legend pos = south west,
			legend cell align=right
		]
		\addplot[red, line width={1pt}, mark=*] table[x index=3,y index=2] {Data/awgm_ExSaw_mw_Run4_conv_info.txt};
		\addlegendentry{AWGM $\norm{\resNE_{k}}$};
		\addplot[red, line width={1pt}, mark=*,dashed] table[x index=3,y index=1] {Data/awgm_ExSaw_mw_Run4_conv_info.txt};
		\addlegendentry{AWGM $\norm{\res_{k}}$};
		\addplot[blue, line width={1pt}, mark=square*] table[x index=3,y index=2] {Data/awgm_ExSaw_SG_mw_conv_info.txt};
		\addlegendentry{SG $\norm{\resNE_{k}}$};
		\addplot[blue, line width={1pt}, mark=square*,dashed] table[x index=3,y index=1] {Data/awgm_ExSaw_SG_mw_conv_info.txt};
		\addlegendentry{SG $\norm{\res_{k}}$};
		\draw (axis cs:6000,0.0005) -- (axis cs:60000,0.00005) -- (axis cs:6000, 0.00005) -- cycle ;
		\node at (axis cs: 12000,0.0001) {{\Large 1}};
		\draw (axis cs:20000,0.02) -- (axis cs:200000,0.02) -- (axis cs:200000, 0.00632455532033676) -- cycle ;
		\node at (axis cs: 110000,0.014) {0.5};
	  \end{axis}
	\end{tikzpicture}\label{Fig:ConvergenceHeat}
	}
	\subfigure[Iteration numbers]{\begin{tikzpicture}[scale=0.52]
	  \begin{axis}
	  		[xlabel={$\mathcal{N}_k = \#\hatbLambda_k$},
			ylabel near ticks,
			legend pos = south east,
			legend cell align=right
		]
		\addplot[red, line width={1pt}, mark=*] table[x index=3,y index=7] {Data/awgm_ExSaw_mw_Run4_conv_info.txt};
		\addlegendentry{AWGM CGLS iterations};
		\addplot[blue, line width={1pt}, mark=square*] table[x index=3,y index=7, skip coords between index={10}{12}] {Data/awgm_ExSaw_SG_mw_conv_info.txt};
		\addlegendentry{SG CGLS iterations};
		\addlegendentry{$1/x$};
	  \end{axis}
	\end{tikzpicture}\label{Fig:ItNbsHeat}
	}
		\subfigure[Test set sizes]{\begin{tikzpicture}[scale=0.52]
	  \begin{axis}
	  		[	xlabel={$\mathcal{N}_k = \#\hatbLambda_k$},
			xmode = log,
			ymode = log,
			ylabel near ticks,
			legend pos = south east,
			legend cell align=right
		]
		\addplot[red, line width={1pt}, mark=*] table[x index=3,y index=4] {Data/awgm_ExSaw_mw_Run4_conv_info.txt};
		\addlegendentry{$\#\checkbLambda_k$};
		\addplot[blue, line width={1pt}, mark=*,dashed] table[x index=3,y index=5] {Data/awgm_ExSaw_mw_Run4_conv_info.txt};
		\addlegendentry{$\#\hatbXi_k$};
		\addplot[blue, line width={1pt}, mark=square*] table[x index=3,y index=6] {Data/awgm_ExSaw_mw_Run4_conv_info.txt};
		\addlegendentry{$\#\checkbXi_k$};

	  \end{axis}
	\end{tikzpicture}\label{Fig:TestSetSizesHeat}
	}
\caption{Heat Equation Example: Comparison of \textbf{LS--AWGM} {(AWGM)} and Sparse Grids {(SG)}.}
\end{figure}

{Starting} with the optimized residual (as in Fig.\ \ref{Fig:XiConstr-Simple}) and the full stable expansions {as in \eqref{Eq:FullStExp}}, we investigate the convergence of the adaptive algorithm and the stability of the finite-dimensional systems \eqref{eq:normalequations_heuristic}. The norms of  primal and  dual residuals are shown in Figure \ref{Fig:ConvergenceHeat} {for AWGM and SG}. As expected, \textbf{LS--AWGM} {reaches} the optimal rate  $s_{\max}=d-1=1$, whereas {uniform SG suffers from the lack of smoothness of the solution.} 
We observe in  Figure \ref{Fig:ItNbsHeat} that the iteration numbers for the least squares {cg}  method in each \textbf{LS--AWGM}-iteration stabilize at about 150 iterations in both approaches. This indicates that the choice of test  sets $\checkbLambda_{k}  = \checkbLambda_{\textrm{Full}}$  yields stability. 
Figure \ref{Fig:TestSetSizesHeat} shows the cardinalities of the test sets. They grow only linearly with $\#\hatbLambda_k$, so that {both $\bw_{\hatbLambda_k}$ and $\resNE_k$} can be computed within linear complexity in each iteration (cf.\ \ref{item:Gal}, \ref{item:Res}). 
\begin{figure}[!ht]
	\subfigure[Convergence of dual residual]{\begin{tikzpicture}[scale=0.52]
	  \begin{axis}
	  		[	xlabel={$\mathcal{N}_k = \#\hatbLambda_k$},
			xmode = log,
			ymode = log,
			ylabel near ticks,
			legend cell align=right
		]
		\addplot[red, line width={1pt}, mark=*] table[x index=3,y index=2] {Data/awgm_ExSaw_mw_Run4_conv_info.txt};
		\addlegendentry{{Optim Constr} $\norm{\resNE_{k}}$};
		\addplot[blue, line width={1pt}, mark=square*] table[x index=3,y index=2] {Data/awgm_ExSaw_mw_Run7_DoubleExp_conv_info.txt};
		\addlegendentry{{Full Constr} $\norm{\resNE_{k}}$};
		\end{axis}
	\end{tikzpicture}\label{Fig:HeatComparisonResConstr-Res}
	}\hspace{0.8cm}
	\subfigure[Index set sizes $\#\hatbXi_{k}$]{\begin{tikzpicture}[scale=0.52]
	  \begin{axis}
	  		[	xlabel={$\mathcal{N}_k = \#\hatbLambda_k$},
			xmode = log,
			ymode = log,
			ylabel near ticks,
			legend pos = south east,
			legend cell align=right
		]
		\addplot[red, line width={1pt}, mark=*] table[x index=3,y index=5] {Data/awgm_ExSaw_mw_Run4_conv_info.txt};
		\addlegendentry{Optim Constr $\#\hatbXi_k$};
		\addplot[blue, line width={1pt}, mark=square*] table[x index=3,y index=5] {Data/awgm_ExSaw_mw_Run7_DoubleExp_conv_info.txt};
		\addlegendentry{Full Constr $\#\hatbXi_k$};
	  \end{axis}
	\end{tikzpicture}\label{Fig:HeatComparisonResConstr-Sizes}
	}
\caption{Heat Equation Example: Comparison of Residual Constructions}\label{Fig:HeatComparisonResConstr}
\end{figure}
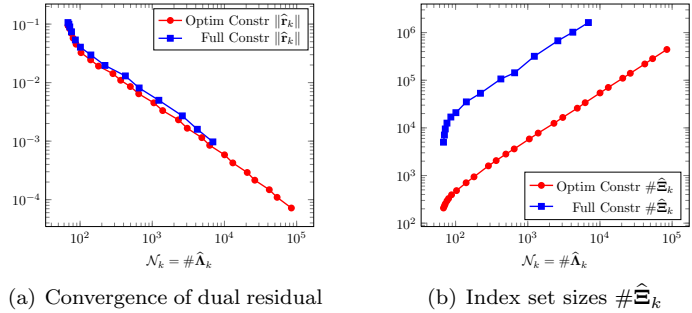
These results are based on {\textbf{OptimResConstruction}} for $\hatbXi_{k}$. In Figure \ref{Fig:HeatComparisonResConstr} {\textbf{FullResConstruction}} is used. As  $\hatbXi_k$ hardly impacts $\res_{k}$, we monitor only the dual residual. {Since}  using a larger index set,  $\norm{\resNE_{k}}_{\ell_2(\hatbcJ)}$ is slightly increased {(as expected)}, but it exhibits the same behaviour as {\textbf{OptimResConstruction}} (Fig.\ \ref{Fig:HeatComparisonResConstr-Res}). This marginal improvement comes at a high cost, $\#\hatbXi_{k}$ is 40-50 times larger, see Fig.\ \ref{Fig:HeatComparisonResConstr-Sizes}. 

\begin{figure}[!ht]
	\subfigure[Residual convergence]{\begin{tikzpicture}[scale=0.52]
	  \begin{axis}
	  		[	xlabel={$\mathcal{N}_k = \#\hatbLambda_k$},
			xmode = log,
			ymode = log,
			ylabel near ticks,
			legend cell align=right,
		]
		\addplot[red, line width={1pt}, mark=*] table[x index=3,y index=2] {Data/awgm_ExSaw_mw_Run4_conv_info.txt};
		\addlegendentry{Full Exp $\norm{\resNE_{k}}$};
		\addplot[red, line width={1pt}, mark=*,dashed] table[x index=3,y index=1] {Data/awgm_ExSaw_mw_Run4_conv_info.txt};
		\addlegendentry{Full Exp $\norm{\res_{k}}$};
		\addplot[blue, line width={1pt}, mark=square*] table[x index=3,y index=2] {Data/awgm_ExSaw_mw_Run5_conv_info.txt};
		\addlegendentry{Red Exp $\norm{\resNE_{k}}$};
		\addplot[blue, line width={1pt}, mark=square*,dashed] table[x index=3,y index=1] {Data/awgm_ExSaw_mw_Run5_conv_info.txt};
		\addlegendentry{Red Exp $\norm{\res_{k}}$};
		\addplot[green, line width={1pt}, mark=diamond*] table[x index=3,y index=2] {Data/awgm_ExSaw_mw_Run6_conv_info.txt};
		\addlegendentry{Temp Exp $\norm{\resNE_{k}}$};
		\addplot[green, line width={1pt}, mark=diamond*,dashed] table[x index=3,y index=1]{Data/awgm_ExSaw_mw_Run6_conv_info.txt};
		\addlegendentry{Temp Exp $\norm{\res_{k}}$};
		\end{axis}
	\end{tikzpicture}\label{Fig:HeatComparisonExpansions-Res}
	}
	\subfigure[Iteration numbers]{\begin{tikzpicture}[scale=0.52]
	  \begin{axis}
	  		[xlabel={$\mathcal{N}_k = \#\hatbLambda_k$},
			ylabel near ticks,
			legend pos = south east,
			legend cell align=right
		]
		\addplot[red, line width={1pt}, mark=*] table[x index=3,y index=7] {Data/awgm_ExSaw_mw_Run4_conv_info.txt};
		\addlegendentry{Full Exp CGLS iterations};
		\addplot[blue, line width={1pt}, mark=square*] table[x index=3,y index=7] {Data/awgm_ExSaw_mw_Run5_conv_info.txt};
		\addlegendentry{Red Exp CGLS iterations};
		\addplot[green, line width={1pt}, mark=diamond*] table[x index=3,y index=7] {Data/awgm_ExSaw_mw_Run6_conv_info.txt};
		\addlegendentry{Temp Exp CGLS iterations};
	  \end{axis}
	\end{tikzpicture}\label{Fig:HeatComparisonExpansions-Its}
	}
	\subfigure[Test set Sszes (normal eq.)]{\begin{tikzpicture}[scale=0.52]
	  \begin{axis}
	  		[	xlabel={$\mathcal{N}_k = \#\hatbLambda_k$},
			xmode = log,
			ymode = log,
			ylabel near ticks,
			legend pos = south east,
			legend cell align=right
		]
		\addplot[red, line width={1pt}, mark=*] table[x index=3,y index=4] {Data/awgm_ExSaw_mw_Run4_conv_info.txt};
		\addlegendentry{Full Exp $\#\checkbLambda_k$};
		\addplot[blue, line width={1pt}, mark=square*] table[x index=3,y index=4] {Data/awgm_ExSaw_mw_Run5_conv_info.txt};
		\addlegendentry{Red Exp $\#\checkbLambda_k$};
		\addplot[green, line width={1pt}, mark=diamond*] table[x index=3,y index=4] {Data/awgm_ExSaw_mw_Run6_conv_info.txt};
		\addlegendentry{Temp Exp $\#\checkbLambda_k$};

	  \end{axis}
	\end{tikzpicture}\label{Fig:HeatComparisonExpansions-CheckLambda}
	}
%
\caption{Heat Equation Example: Comparison of Stable Expansions}\label{Fig:HeatComparisonExpansions}
\end{figure}
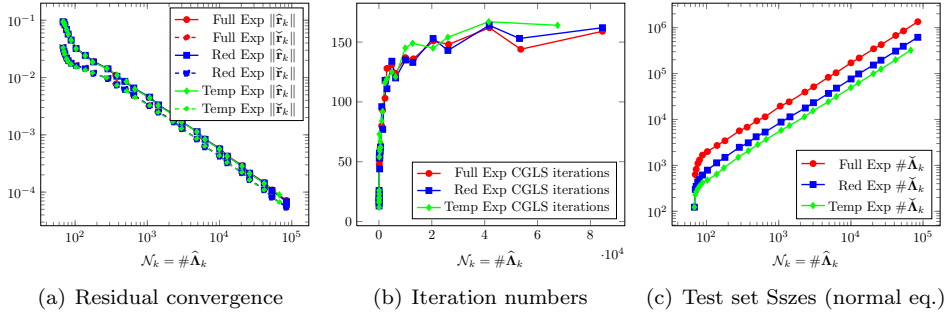

Finally, in Figure \ref{Fig:HeatComparisonExpansions}, we compare the stable expansion types {(Full, Reduced, Temporal)}. {We} find no discernible differences in the residual  (Fig.\ \ref{Fig:HeatComparisonExpansions-Res}) and only a very slight increase in the iteration numbers in  \textbf{GALSOLVE}  (Fig.\ \ref{Fig:HeatComparisonExpansions-Its}).
It seems that choosing $\checkbXi_k = \textbf{TemporalStableExpansion}(\hatbXi^\text{tmp}_k,1)$ yields results that are comparable to the other extensions, which  could not be deduced from \cite{Kestler:2012c,Kestler:2012d}.
All three methods seem stable, and we can  reduce the size of the test  sets by a factor of about 3.4 for $\checkbLambda_k$ (and likewise by 2.5 for $\checkbXi_k$).

\subsection{Convection-Diffusion-Reaction Equation}
\begin{figure}[!hb]
	\subfigure[Solution $u(t,x)$ of the CDR example]{
	\includegraphics[width=0.4\textwidth]{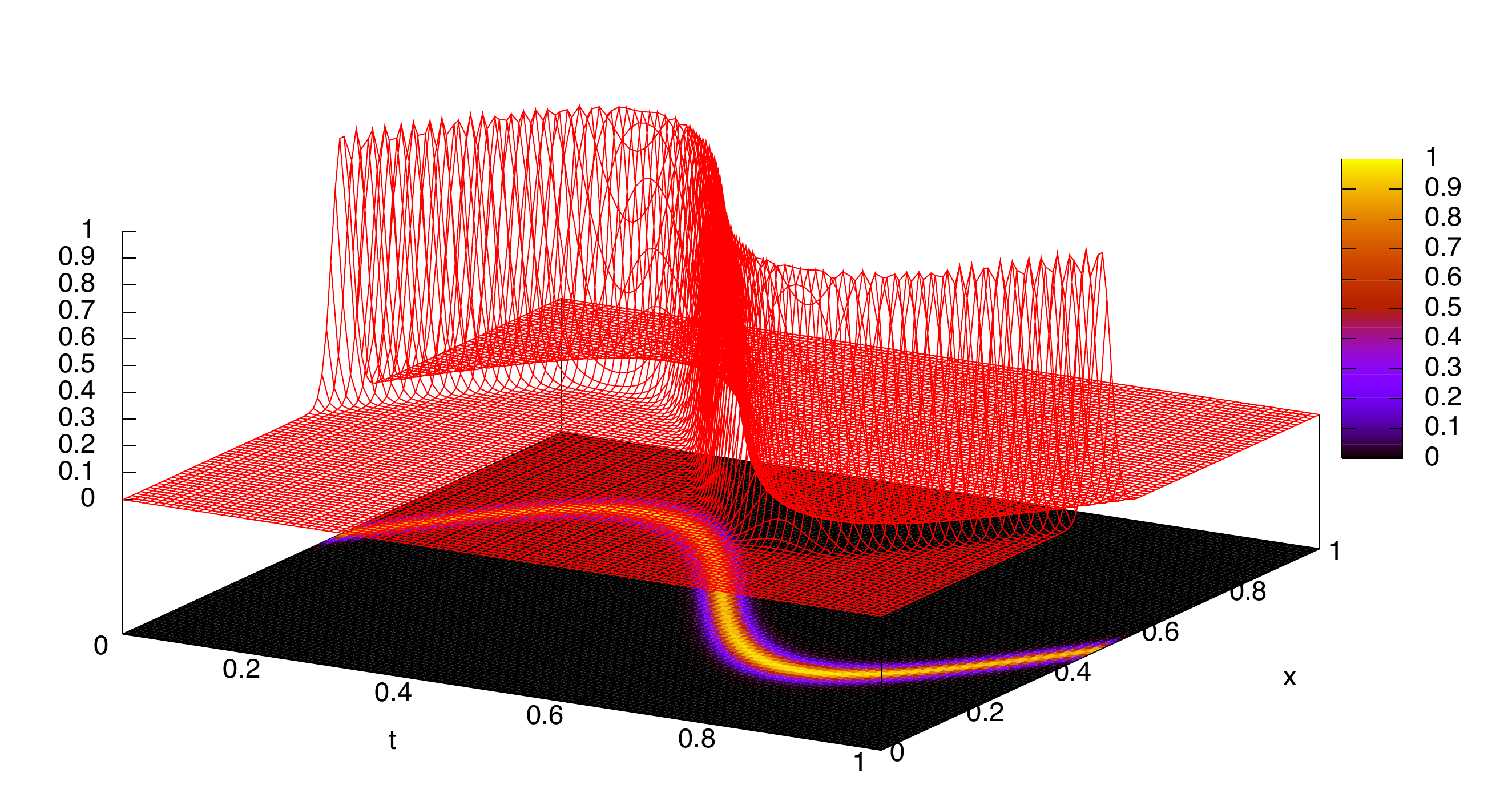}
		\label{Fig:ExCDRCurve-Sol}
 	}\hspace{0.8cm}
	\subfigure[Support centres of basis functions in $\hatbLambda_{k}$ for $k=12$ ($\#\hatbLambda_{12}=9445$).]{
		\includegraphics[width=0.4\textwidth]{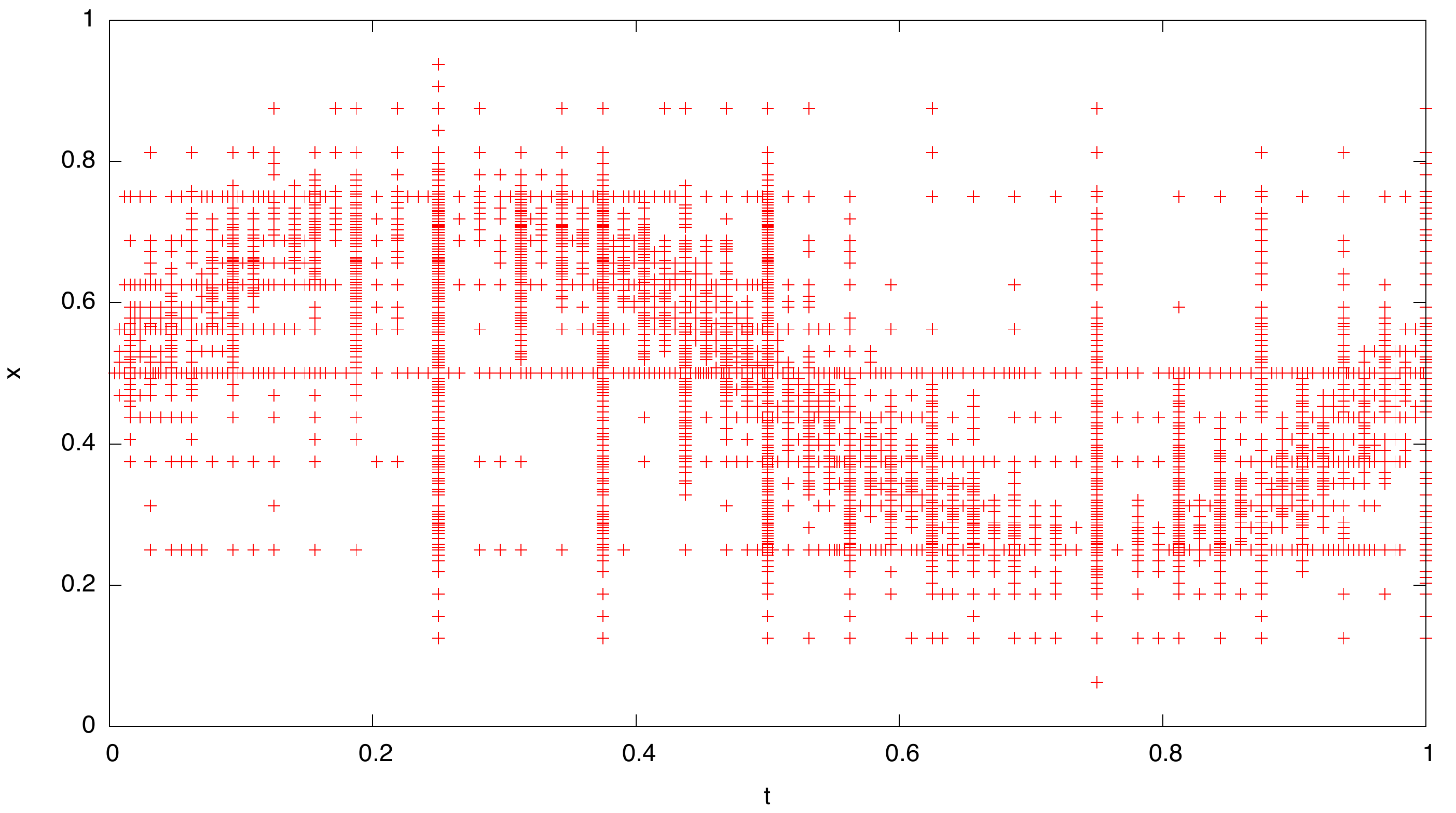}
	\label{Fig:ExCDRCurve-SupportCentres}
	}
\caption{CDR Example: Solution and adaptive refinement}
\end{figure}

As a second example, we {consider} the convection-diffusion-reaction (CDR) equation 
\begin{equation*}
	\left\{
	\begin{aligned}
	u_{t} - u_{xx} + u_{x} + u &= f(t,x) \quad\qquad \text{ on } \Omega = (0,1), \\[1ex]
			u(t,0) &= u(t,1) \hspace{1.1cm} \text{ for all } t \in [0,T],\\
			u(0,x) = u(T,x)&= 0 \hspace{2.05cm} \!\text{on } \overline\Omega,
	\end{aligned}\right.
\end{equation*}
for a $f(t,x)$ that yields $u(t,x) = e^{-1000 \left(x - (0.5 + 0.25 \sin(2 \pi t))\right)^{2}}$, see Figure \ref{Fig:ExCDRCurve-Sol}.
Note that $u$ is infinitely smooth but exhibits large gradients in non axis-aligned directions.

The support centers (i.e., the centers of $\supp \hatbpsi_{\blambda}$, $\blambda \in \hatbLambda_k$) in Figure \ref{Fig:ExCDRCurve-SupportCentres} indicate that  {the AWGM benefits} from its ability to refine not only independently in each dimension, but in particular \emph{locally in the full space-time domain}. 
\begin{figure}[!ht]
	\subfigure[Residual convergence]{
	\begin{tikzpicture}[scale=0.52]
	  \begin{axis}
	  		[	xlabel={$\mathcal{N}_k = \#\hatbLambda_k$},
			xmode = log,
			ymode = log,
			ylabel near ticks,
			legend pos = south west,
			legend cell align=left
			]
			\addplot[red, line width={1pt}, mark=*] table[x index=3,y index=2] {Data/awgm_ExCDRCurve_mw_Run5_Temporal_conv_info.txt};
	   	\addlegendentry{AWGM $\norm{\resNE_{k}}$};
			\addplot[red, line width={1pt}, mark=*,dashed] table[x index=3,y index=1] {Data/awgm_ExCDRCurve_mw_Run5_Temporal_conv_info.txt};
		\addlegendentry{AWGM $\norm{\res_{k}}$};
			\addplot[blue, line width={1pt}, mark=square*] table[x index=3,y index=2] {Data/awgm_ExCDRCurve_SG_mw_conv_info.txt};
		  \addlegendentry{SG $\norm{\resNE_{k}}$};
		  \addplot[blue, line width={1pt}, mark=square*,dashed] table[x index=3,y index=1] {Data/awgm_ExCDRCurve_SG_mw_conv_info.txt};
		  \addlegendentry{SG $\norm{\res_{k}}$};
			\draw (axis cs:5000,0.3) -- (axis cs:50000,0.03) -- (axis cs:5000, 0.03) -- cycle ;
			\node at (axis cs: 10000,0.06) {1};
		\end{axis}
	\end{tikzpicture}\label{Fig:ExCDRCurve-Res}
	}\hspace{0.8cm}
	\subfigure[Iteration numbers]{\begin{tikzpicture}[scale=0.52]
	  \begin{axis}
	  		[xlabel={$\mathcal{N}_k = \#\hatbLambda_k$},
			ylabel near ticks,
			legend pos = south east,
			legend cell align=left
		]
		\addplot[red, line width={1pt}, mark=*] table[x index=3,y index=7] {Data/awgm_ExCDRCurve_mw_Run5_Temporal_conv_info.txt};
		\addlegendentry{AWGM CGLS iterations};
	  \end{axis}
	\end{tikzpicture}\label{Fig:ExCDRCurve-Its}
	}
	\caption{CDR Example: {Convergence and stability of \textbf{LS-AWGM} (AWGM) and Sparse Grids (SG)}}
\end{figure}
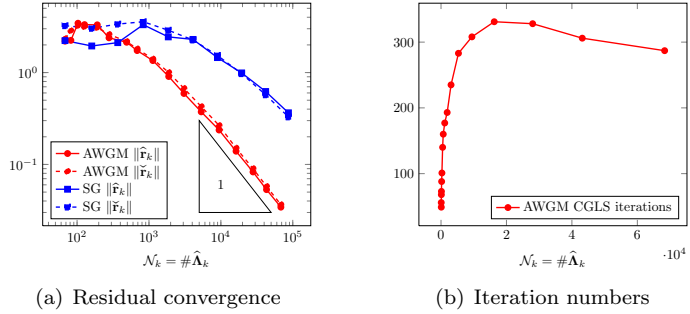
This is also mirrored in Figure \ref{Fig:ExCDRCurve-Res}, where we observe  the optimal $s_{\max}=1$ for the \textbf{LS--AWGM}, and a stable number of inner iterations (Fig.~\ref{Fig:ExCDRCurve-Its}), employing the optimized construction of $\hatbXi_{k}$ and only temporal stable expansions for $\checkbLambda_{k}$, $\checkbXi_{k}$. The smoothness of the solution allows for a convergence rate close to 1 for the sparse grid approach, however, the asymptotic regime and comparable residual norms are only reached for index sets that are over a magnitude larger.

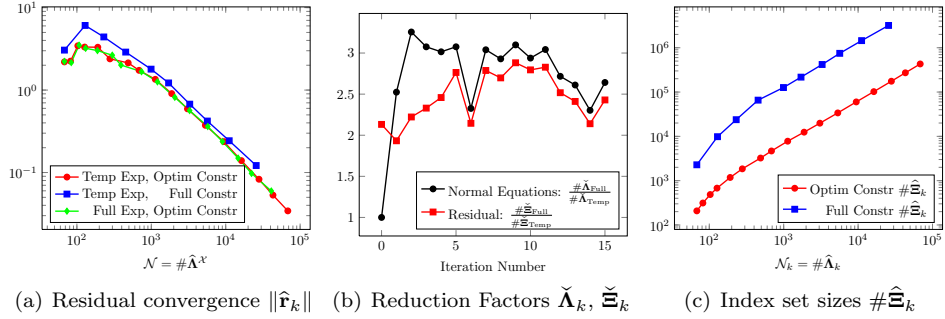
\begin{figure}[!ht]
	\subfigure[Residual convergence $\norm{\resNE_{k}}$]{\begin{tikzpicture}[scale=0.52]
	  \begin{axis}
	  		[	xlabel={$\mathcal{N} = \#\hatbLambda^{\cX}$},
			xmode = log,
			ymode = log,
			ylabel near ticks,
			legend pos=south west,
			legend cell align=right
		]
		\addplot[red, line width={1pt}, mark=*] table[x index=3,y index=2] {Data/awgm_ExCDRCurve_mw_Run5_Temporal_conv_info.txt};
		\addlegendentry{Temp Exp, Optim Constr};
		\addplot[blue, line width={1pt}, mark=square*] table[x index=3,y index=2] {Data/awgm_ExCDRCurve_mw_Run6_DoubleExp_conv_info.txt};
		\addlegendentry{Temp Exp, $\quad\, $Full Constr};
		\addplot[green, line width={1pt}, mark=diamond*] table[x index=3,y index=2] {Data/awgm_ExCDRCurve_mw_Run4_conv_info.txt};
		\addlegendentry{Full Exp, Optim Constr};

		\end{axis}
	\end{tikzpicture}\label{Fig:CDRComparisonExpansions-Res}
	}
	\subfigure[Reduction Factors $\checkbLambda_{k}$, $\checkbXi_{k}$ ]{\begin{tikzpicture}[scale=0.52]
	  \begin{axis}
	  		[	xlabel={Iteration Number},
			ylabel near ticks,
			legend pos = south east,
			legend cell align=left,
			legend style={inner ysep=4pt}
		]
		\addplot[black, line width={1pt}, mark=*] table[x index=0,y expr=\thisrowno{4}/\thisrowno{12}] {Data/awgm_ExCDR_conv_info_Full_Temporal.txt};
		\addlegendentry{Normal Equations: $\tfrac{\#\checkbLambda_{\textrm{Full}}}{\#\checkbLambda_{\textrm{Temp}}}$};
		\addplot[red, line width={1pt}, mark=square*] table[x index=0,y expr=\thisrowno{6}/\thisrowno{14}] {Data/awgm_ExCDR_conv_info_Full_Temporal.txt};
		\addlegendentry{Residual: $\tfrac{\#\checkbXi_{\textrm{Full}}}{\#\checkbXi_{\textrm{Temp}\vspace{2ex}}}$};
	  \end{axis}
	\end{tikzpicture}\label{Fig:CDRComparisonExpansions-RedFactors}\vspace{0.1cm}
	}
		\subfigure[Index set sizes $\#\hatbXi_{k}$]{\begin{tikzpicture}[scale=0.52]
	  \begin{axis}
	  		[	xlabel={$\mathcal{N}_k = \#\hatbLambda_k$},
			xmode = log,
			ymode = log,
			ylabel near ticks,
			legend pos = south east,
			legend cell align=right
		]
		\addplot[red, line width={1pt}, mark=*] table[x index=3,y index=5] {Data/awgm_ExCDRCurve_mw_Run5_Temporal_conv_info.txt};
		\addlegendentry{Optim Constr $\#\hatbXi_k$};
		\addplot[blue, line width={1pt}, mark=square*] table[x index=3,y index=5] {Data/awgm_ExCDRCurve_mw_Run6_DoubleExp_conv_info.txt};
		\addlegendentry{Full Constr $\#\hatbXi_k$};
	  \end{axis}
	\end{tikzpicture}\label{Fig:CDRComparisonConstr-HatXi}
	}
\caption{CDR Example: Comparison of different index set constructions.}\label{Fig:CDRComparisonExpansions}
\end{figure}

Finally, we compare the above AWGM results with those obtained for larger sets, i.e., using full stable expansions and the \textbf{FullResConstruction} for $\hatbXi_{k}$.
As {before}, we see in Figures ~\ref{Fig:CDRComparisonExpansions-Res}, \ref{Fig:CDRComparisonExpansions-RedFactors} that we can reduce the size of the test  sets $\checkbLambda_k$, $\checkbXi_k$ by factors  2 to 3 without losing accuracy. Likewise, the full construction of $\hatbXi_{k}$ yields index sets that are approximately 20 times as large as for the optimized version (cf. Fig.~\ref{Fig:CDRComparisonConstr-HatXi}), with only a slight improvement in the residual approximation.

\appendix
\section{{Proof of Proposition \ref{Prop:2.1}}}\label{sec:appendix-wellposed}

{
\noindent We follow \cite{Schwab:2009} to verify the Babu\v{s}ka-Aziz conditions in a time-periodic setting.

\paragraph*{\emph{(1) Continuity.}}
Follows from \eqref{eq:ellipticityA}, the definitions of $\norm{\cdot}_{\cX}$, $\norm{\cdot}_{\cY}$ as well as Cauchy-Schwarz's, H\"older's and Young's inequalities.

\paragraph*{\emph{(2) Inf-sup condition.}}
We consider an arbitrary $0\ne w \in \cX$ and define $z_w(t) := \left(\OpSp^*\right)^{-1} \dot{w}(t)$ for the adjoint $\OpSp^{*}$ of $\OpSp$. 
The bound $\norm{(\OpSp^{*})^{-1}} \leq \alpha^{-1}$ then yields for $v_w(t) := z_w(t) + w(t)$ that $\norm{v_w}_{\cY}\leq\ \sqrt{2} \max\{1, \alpha^{-1} \} \norm{w}_{\cX} < \infty$. 
By definition of $z_{w}$ and \eqref{eq:ellipticityA}, $\eval{z_w(t)}{\dot{w}(t)}_{V\times V'}= \eval{z_w(t)}{\OpSp[z_w(t)]}_{V\times V'} \geq \alpha \norm{z_w(t)}_V^2 \geq \frac{\alpha}{\gamma^2} \norm{\dot{w}(t)}_{V'}^2$. Since $w \in \cX$ is periodic, we have  $\int_{0}^{T}\eval{w}{\dot w}_{V\times V'} + \eval{z_{w}}{\OpSp[w]}_{V\times V'} dt= \int_{0}^{T}\eval{w}{\dot w}_{V\times V'} + \int_{0}^{T}\eval{\dot{w}}{w}_{V\times V'}= \int_{0}^{T} \frac{d}{dt}\norm{w(t)}_{H}^{2}dt \\= \norm{w(T)}_{H}^{2}-\norm{w(0)}_{H}^{2} = 0$, so that we finally get 
$b(w, v_w) 
 \ge \alpha \min\{1, \gamma^{-2}\}\norm{w}_{\cX}^{2} \ge \frac{\alpha \min\{1, \gamma^{-2}\} }{\sqrt{2}  \max\{1, \alpha^{-1}\}} \norm{w}_{\cX}\norm{v_w}_{\cY} > 0$.

\paragraph*{\emph{(3) Surjectivity}}

Let $0\ne v \in \cY$. We aim to construct $z \in \cX$ with $\eval{w(t)}{\dot z(t)}_{V\times V'} + \eval{w(t)}{\OpSp[z(t)]}_{V\times V'} = \eval{w(t)}{\OpSp[v(t)]}_{V\times V'}$ for all $w \in \cY$, and $t$ a.e.\ on $(0,T)$, as then  $b(z,v) = \int_{0}^{T}\eval{v(t)}{\OpSp[v(t)]}_{V\times V'} \geq \alpha \norm{v}_{\cY}^{2} > 0$,  so that the surjectivity condition is fulfilled.

\subparagraph{\it (i) Faedo-Galerkin approximation of an initial value problem}
Let $\{\phi_i:i \in \N\}$ be a basis for $V$, $V_n := \Span\{\phi_i, i=1,\dots,n\}$, $z_n(t) := \sum_{i=1}^n z_i^{(n)}(t) \phi_i$. Then the linear system of ODEs
$\eval{ w_n}{\dot z_n(t)}_{V\times V'} + \eval{w_{n}}{\OpSp[z_{n}(t)]}_{V\times V'} = \eval{w_{n}}{\OpSp[v(t)]}_{V\times V'}$, $z_n(0) = z_{n0}$, 
has a solution  $z_n \in C(0,T; V_n)$ with  $\dot z_n \in L_2(0,T;V_n)$ for all $w_n \in V_n$  a.e.\ on $I$ and for (arbitrary) $z_0 \in H$ and its orthogonal projection $z_{n0}$ onto $V_{n}$. 

\subparagraph{\it (ii)  A-priori estimates} 
(i), \eqref{eq:ellipticityA} and Young's inequality with some $\eps < \frac{\alpha}{\gamma}$ yield
\begin{align}
\frac12 \frac{d}{dt} \norm{z_n(t)}_{H}^{2}+ \alpha \norm{z_n(t)}_{V}^{2}
	&\leq \gamma \eps\norm{z_n(t)}_V^2 + \frac{\gamma}{4\eps}\norm{v(t)}_V^2  \label{eq:BoundAPriori}
 \end{align}
and hence by integration over $[0,s]$, $s \in [0,T]$, using $(\alpha - \gamma \eps)>0$ that
$
	\norm{z_{n}(s)}_{H}^{2} - 
	\norm{z_{n}(0)}_{H}^{2} = \int_{0}^{s}\frac{d}{dt} \norm{z_n(t)}_{H}^{2} dt
	\ \leq \ \frac{\gamma}{2\eps} \int_{0}^{s}\norm{v(t)}_{V}^{2} dt, 
$
so that $\sup_{s \in [0,T]} \norm{z_{n}(s)}_{H}^{2} < \infty$ and $\{z_{n}\}_{n\in \N}$ is uniformly bounded in $L_{\infty}(0,T;H)$. Similarly, we can conclude that $2(\alpha - \gamma\eps) \| z_n\|_{\cY}  \leq  \norm{z_{n}(0)}_{H}^{2} - \norm{z_{n}(T)}_{H}^{2} + \frac{\gamma}{2\eps} \norm{v}_{\cY}^{2} <\ \infty$,
so that $\{z_{n}\}_{n\in \N}$ is also uniformly bounded in $\cY$.

\subparagraph{\it (iii) Periodicity} 
Abbreviate $\bar c := \frac{\gamma}{4\eps}$, $\bar \alpha := 2\frac{(\alpha - \gamma \eps)}{c_{1}} >0$ with $ c_{1} := \sup_{\phi \in V} \frac{\norm{\phi}_{V}}{\norm{\phi}_{H}}$ and
multiply \eqref{eq:BoundAPriori} by $e^{\bar\alpha t}$. Then 
$
\frac{d}{dt} \left(e^{\bar\alpha t}\,\norm{z_{n}(t)}_{H}^{2} \right)
	= e^{\bar\alpha t}\frac{d}{dt} \norm{z_{n}(t)}_{H}^{2} + e^{\bar\alpha t} \bar \alpha\, \norm{z_{n}(t)}_{H}^{2}
	\le e^{\bar\alpha t} \bar c\, \norm{v(t)}_{V}^{2}$  
and by integration over $[0,T]$, we obtain 
\begin{equation}\label{Eq:BoundDiscreteSolAtFinalTime}
	\norm{z_{n}(T)}_{H}^{2} \ 
	\leq \ e^{- \bar\alpha T}\norm{z_{n}(0 )}_{H}^{2}+ \bar c\, e^{-\bar\alpha T}\int_{0}^{T} e^{\bar\alpha t}\norm{v(t)}_{V}^{2}\, dt.
\end{equation}
Set $M:= \{ z \in V_{n} : \norm{z}_{H} \leq R := K^{\frac 12} (1 - e^{- \bar\alpha T})^{-\frac 1 2}\}$, $K := \bar ce^{-\bar\alpha T}\int_{0}^{T} e^{\bar\alpha t}\norm{v(t)}_{V}^{2}dt$.
The set $M$ is convex and compact in $V_{N}$. If $z_{n}(0) \in M$, \eqref{Eq:BoundDiscreteSolAtFinalTime} implies that $\norm{z_{n}(T)}_{H}^{2} \leq e^{- \bar\alpha T} R^{2} + K \leq R$, i.e. $z_{n}(T) \in M$.
 Since by Gronwall's lemma the mapping $S: M \to M$, $z_{n}(0) \mapsto z_{n}(T)$, is continuous, the existence of a fixed-point $S(\bar z_{n}) = \bar z_{n} \in M$ follows from Brouwer's fixed-point theorem. By the a-priori estimates, the sequence $\{\bar z_{n}\}_{n\in \N}$ is bounded in $H$, so that there exists a subsequence (also denoted by $\{\bar z_{n}\}$) converging weakly to some $\bar z \in H.$
\smallskip
%

\paragraph{\it (iv) Convergence} 
Consider the periodic solution $z_{n}(t)$ from (iii), i.e. the solution of the ODE system with initial value $z_{n0} = \bar z_{n}$. From the a-priori estimates, we have that $\{z_{n}\}$ is uniformly bounded in the separable space $\cY$, so that there exists a subsequence (also denoted $\{z_{n}\}$) converging weakly to some $z$ in $\cY$. 
For $w_{n} := \theta(t)\phi_{j}$, $\theta(t) \in C^{1}(0,T)$, we then have by integration over $[0,T]$ and integration by parts of the first term that for all $j=1,\dots,n$
$
-\eval{\dot\theta \phi_{j}}{z_{n}}
 	= \eval{\theta(0)\phi_{j} - \theta(T)\phi_{j}}{\bar z_{n}}_{H} + \eval{\theta\phi_j}{\OpSp[v-z_n]}
$.
As $z_{n} \rightharpoonup z$ in $\cY$ and $\bar z_{n} \rightharpoonup \bar z$ in $H$, we can pass to the limit $n \to \infty$ and obtain
\begin{equation}\label{eq:partInt}
-\eval{\dot\theta \phi_{j}}{z}
 	= \eval{\theta(0)\phi_{j} - \theta(T)\phi_{j}}{\bar z}_{H} + \eval{\theta\phi_j}{\OpSp[v-z]}.
\end{equation}
This particularly holds true for all $\theta \in \mathcal{D}(I)$, so that $\dot z = \mathcal{A}(\cdot)(v-z)$ in the distributional sense and hence $\dot z \in L_{2}(0,T;V')$. 
Moreover, (\ref{eq:partInt}) implies that for $w \in C^{1}(0,T;V)$, we have 
$-\eval{ \dot w}{z} - \eval{w(0) - w(T)}{\bar z} = \eval{w}{\OpSp[v-z]}  = \eval{ \dot z}{w} = -\eval{ \dot w}{z} +  \eval{w(T)}{z(T)}_{H}  - \eval{w(0)}{z(0)}_{H}$, 
so that indeed $\bar z = z(0) = z(T)$ in $H$ and hence $z \in \cX$. With this $z$, the surjectivity condition is fulfilled.

}

\section{Auxiliary wavelet compression results}\label{sec:appendix}

Here, we {report two}  facts for $\bB$ defined in \eqref{eq:equiv_ell2_prob} which are required in Section \ref{sec:impl}. {We shall always assume}  that \eqref{eq:compr} {holds. For further details, we refer to \cite{Kestler:Diss}.} 
\begin{lemma}[\cite{Kestler:2012d}] \label{lem:comr_bBtrbB}
For sufficiently small $\eta<1$,  $\bB_\eta \in \cL(\ell_2(\hatbcJ), \ell_2(\checkbcJ))$ and $ \bB^\top_\eta \bB_\eta \in \cL( \ell_2(\hatbcJ) , \ell_2(\hatbcJ))$ are boundedly invertible with bounds depending on $\eta$.
\end{lemma}
%


\begin{proposition}[\cite{Kestler:2012d}] \label{prop:dual_res}
Let the assumptions of Theorem \ref{thm:primal_res} hold.
Then, there exists a constant $\omegals$ such that
$\| \bB^\top(\bf - \bB \bw_{\hatbLambda}) - \resNE \|_{\ell_2(\hatbcJ)} \leq \omegals \| \resNE \|_{\ell_2(\hatbcJ)}$, $\resNE := \bB_\eta^\top \res$.
\end{proposition}



\bibliographystyle{alpha}
\bibliography{literature}

\end{document}

%% file: newcommands.tex
\DeclareFontFamily{U}{mathx}{\hyphenchar\font45}
\DeclareFontShape{U}{mathx}{m}{n}{
     <5> <6> <7> <8> <9> <10>
      <10.95> <12> <14.4> <17.28> <20.74> <24.88>
      mathx10
      }{}
\DeclareSymbolFont{mathx}{U}{mathx}{m}{n}
\DeclareFontSubstitution{U}{mathx}{m}{n}
\DeclareMathAccent{\mywidehat}{0}{mathx}{"70} 
\DeclareMathAccent{\mywidecheck}{0}{mathx}{"71}

\newcommand{\argmin}{\operatornamewithlimits{argmin}}

\DeclareMathOperator{\diam}{diam}
\DeclareMathOperator{\diag}{diag}
\DeclareMathOperator{\supp}{supp}
\DeclareMathOperator{\Span}{span}
\DeclareMathOperator{\dist}{dist}

\newcommand{\triplenorm}{|\!|\!|}
\newcommand{\cAmtree}{\mathcal{A}_{\mathrm{mtree}}}
\newcommand{\omegals}{\omega_{\mathrm{ls}}}
\newcommand{\gammals}{\gamma_{\mathrm{ls}}}
\newcommand{\sstar}{s^\ast}

\newcommand{\per}{\mathrm{per}}
\newcommand{\Id}{\mathrm{Id}}
\newcommand{\smax}{s_{\max}}

\newcommand{\scrH}{\mathscr{H}}

\newcommand{\R}{\mathbb{R}}
\newcommand{\N}{\mathbb{N}}

\newcommand{\cA}{\mathcal{A}}
\newcommand{\cB}{\mathcal{B}}
\newcommand{\cC}{\mathcal{C}}

\newcommand{\cH}{\mathcal{H}}

\newcommand{\cJ}{\mathcal{J}}

\newcommand{\cL}{\mathcal{L}}
\newcommand{\cN}{\mathcal{N}}

\newcommand{\cO}{\mathcal{O}}

\newcommand{\cX}{\mathcal{X}}
\newcommand{\cY}{\mathcal{Y}}

\renewcommand{\bf}{\mathbf{f}}
\newcommand{\bg}{\mathbf{g}}
\newcommand{\br}{\mathbf{r}}
\newcommand{\bu}{\mathbf{u}}
\newcommand{\bv}{\mathbf{v}}
\newcommand{\bw}{\mathbf{w}}
\newcommand{\bx}{\mathbf{x}}

\newcommand{\bB}{\mathbf{B}}
\newcommand{\bC}{\mathbf{C}}
\newcommand{\bD}{\mathbf{D}}
\newcommand{\bE}{\mathbf{E}}
\newcommand{\bR}{\mathbf{R}}

\newcommand{\bcJ}{\bm{\cJ}}
\newcommand{\blambda}{\boldsymbol \lambda}
\newcommand{\bLambda}{\boldsymbol \Lambda}
\newcommand{\bmu}{\boldsymbol \mu}
\newcommand{\bOmega}{\boldsymbol \Omega}
\newcommand{\bpsi}{\boldsymbol \psi}
\newcommand{\bPsi}{\boldsymbol \Psi}
\newcommand{\bPi}{\boldsymbol \Pi}
\newcommand{\bsigma}{\boldsymbol \sigma}
\newcommand{\bSigma}{\boldsymbol \Sigma}
\newcommand{\bXi}{\boldsymbol \Xi}

\newcommand{\bcJx}{\bcJ_{\hspace{-0.5ex}x}}
\newcommand{\hatbcJ}{\hspace{0.9ex}\mywidehat{\hspace{-0.9ex}\bm{\cJ}}}
\newcommand{\checkbcJ}{\hspace{0.9ex}\mywidecheck{\hspace{-0.9ex}\bm{\cJ}}}

\newcommand{\checkbpsi}{\mywidecheck{\bpsi}}
\newcommand{\checkpsi}{\mywidecheck{\psi}}
\newcommand{\hatbpsi}{\mywidehat{\bpsi}}
\newcommand{\hatpsi}{\mywidehat{\psi}}
\newcommand{\checkPsi}{\mywidecheck{\Psi}}
\newcommand{\hatPsi}{\mywidehat{\Psi}}
\newcommand{\checkbPsi}{\mywidecheck{\bPsi}}
\newcommand{\hatbPsi}{\mywidehat{\bPsi}}

\newcommand{\checkbLambda}{\mywidecheck{\bLambda}}
\newcommand{\hatbLambda}{\mywidehat{\bLambda}}

\newcommand{\checkLambda}{\mywidecheck{\Lambda}}
\newcommand{\hatLambda}{\mywidehat{\Lambda}}

\newcommand{\checkbXi}{\mywidecheck{\bXi}}
\newcommand{\hatbXi}{\mywidehat{\bXi}}

\newcommand{\resNE}{\mywidehat{\br}}
\newcommand{\res}{\mywidecheck{\br}}

\newcommand{\eps}{\varepsilon}
\newcommand{\wt}{\widetilde}

\newcommand{\vecA}{\vec{A}}
\newcommand{\vecB}{\vec{B}}
\newcommand{\vecId}{\vec{\Id}}
\newcommand{\vecL}{\vec{L}}
\newcommand{\vecU}{\vec{U}}

\newcommand{\eval}[2]{\langle #1 , #2 \rangle}
\newcommand{\bbB}[2]{{}_{#1} \bB_{#2}}
\newcommand{\bbBB}[2]{{}_{#1} \bB^\top \bB_{#2}}

\newcommand{\bbBtr}[2]{{}_{#1} \bB_{#2}{}^{\hspace{-1.8ex}\top}}
\newcommand{\rmc}[2]{\mathrm{c}_{#1}(#2)}
\newcommand{\rmC}[2]{\mathrm{C}_{#1}(#2)}
\newcommand{\rmD}[1]{\mathrm{D}_{#1}}
\newcommand{\norm}[1]{\lVert #1\rVert}

\newcommand{\Domain}{\bOmega}
\newcommand{\OpSp}{\cA(t)}
\newcommand{\OpST}{\cB}
\newcommand{\gammaB}{\gamma_{\OpST}} 
\newcommand{\betaB}{\beta_{\OpST}} 